\documentclass[12pt]{amsart}
\usepackage{amsmath,amssymb, xypic,
verbatim,amscd}
\usepackage[dvips]{graphics}
\numberwithin{equation}{section}

\newcommand\temp{t}
\newcommand\Ftwo{\mathcal{K}}
\newcommand\SL{\operatorname{SL}}
\newcommand\Tor{\operatorname{Tor}}

\newcommand\cee{\chi}
\newcommand\ess{C}
\newcommand\dm{\mathfrak{d}}
\newcommand\Quot{\operatorname{Quot}}
\newcommand\yess{S}
\newcommand\ttt{(i^p_{\alpha^p(t)}-\alpha^p(t))}

\newcommand\mfm{\operatorname{Gr}}

\newcommand\nation{\delta}
\newcommand\mupar{\mu_{\operatorname{par}}}
\newcommand\SU{\operatorname{SU}}
\newcommand\Unit{\operatorname{U}}

\newcommand\shom{\mathcal{H}om}
\newcommand\send{\mathcal{E}{nd}}

\newcommand\pa{\operatorname{pardeg}}
\newcommand\my{\mathcal{Y}}

\newcommand\mug{\muu_{\operatorname{gen}}(\mathcal{V},\mathcal{Q},\mathcal{I})}

\newcommand\Disc{\operatorname{Disc}}

\newcommand\Hone{\operatorname{H}^1}
\newcommand\ext{\operatorname{Ext}}
\newcommand\homo{\operatorname{Hom}}
\newcommand\Hom{\operatorname{Hom}}
\newcommand\tensor{\otimes}
\newcommand\im{\operatorname{im}}
\newcommand\Ext{\operatorname{Ext}}

\newcommand\codim{{\operatorname{codim}}}
\newcommand\dimm{{\dim}}
\newcommand\Fl{\operatorname{Fl}_{\mpp}}
\newcommand\Flag{\operatorname{Fl}}
\newcommand\Gl{\operatorname{GL}}
\newcommand\Spec{\operatorname{Spec}}
\newcommand\rk{\operatorname{rk}}

\newcommand\umw{\underline{\mw}}
\newcommand\wt{\operatorname{wt}}
\newcommand\tmw{\tilde{\mathcal{W}}}
\newcommand\tmq{\tilde{\mathcal{Q}}}
\newcommand\tr{\tilde{r}}
\newcommand\td{\tilde{d}}
\newcommand{\leto}[1]{\stackrel{#1}{\to}}
\newcommand{\letoo}[1]{\stackrel{#1}{\hookrightarrow}}

\newcommand\Gr{\operatorname{Gr}}

\newcommand\ma{\mathcal{A}}
\newcommand\mb{\mathcal{B}}
\newcommand\mc{\mathcal{C}}

\newcommand\me{\mathcal{E}}
\newcommand\mee{\mathcal{M}}
\newcommand\mh{\mathcal{H}}

\newcommand\mf{\mathcal{F}}
\newcommand\mg{\mathcal{G}}
\newcommand\mz{\mathcal{Z}}
\newcommand\mi{\mathcal{I}}
\newcommand\mj{\mathcal{J}}
\newcommand \braI{\langle\mi\rangle}
\newcommand\mk{\mathcal{K}}
\newcommand\ml{\mathcal{L}}
\newcommand\mpp{{\operatorname{ S}}}
\newcommand\mq{\mathcal{Q}}

\newcommand\ms{\mathcal{S}}
\newcommand\mt{\mathcal{T}}

\newcommand\muu{{\Bbb{U}}}
\newcommand\mv{\mathcal{V}}
\newcommand\tmv{\tilde{\mathcal{V}}}
\newcommand\umv{\underline{\mathcal{V}}}
\newcommand\mw{\mathcal{W}}
\newcommand\hoo{\Hom_{\mi}(\mv,\mq,\mf,\mg)}
\newcommand\hoos{\Hom_{\mk}(\ms,\mv/\ms,\mf(\ms),\mf(\mv/\ms))}
\newcommand\hos{T}
\newcommand\cn{\Bbb{C}^n}
\newcommand\sS{{\operatorname{Sh} }}

\newcommand\bull{\sssize{\bullet}}
\newcommand\pone{\operatorname{\Bbb{P}^1}}

\newcommand\dimii{\dim \mi }
\newcommand\dimkk{\dim \mk }
\newcommand\bed{\beta}
\newcommand\emm{\Lambda  }
\newcommand\bedda{\beta}
\newcommand\emman{\Lambda  }
\newcommand\bedd{D-d}
\newcommand\emma{n-r  }
\newcommand\degr{\delta}
\newcommand\hoq{\Hom(\mv,\mq,\td,\tr,\epsilon)}
\newcommand\gammar{\gamma }
\newtheorem{theorem}{Theorem}[section]
\newtheorem{lemma}[theorem]{Lemma}

\newtheorem{proposition}[theorem]{Proposition}
\newtheorem{corollary}[theorem]{Corollary}
\newtheorem{definition-proposition}[theorem]{Definition-Proposition}

\theoremstyle{definition} \newtheorem{definition}[theorem]{Definition}
\theoremstyle{defi} \newtheorem{defi}[theorem]{Definition}

 \newtheorem{claim}[theorem]{Claim}

\theoremstyle{remark} \newtheorem{remark}[theorem]{Remark}
\begin{document}
\title[The Quantum Horn Problem]{Quantum Generalization of the Horn Conjecture}
\author{Prakash Belkale}
\address{Department of Mathematics\\ UNC-Chapel Hill\\ CB \#3250, Phillips Hall
\\ Chapel Hill , NC 27599}
\email{belkale@email.unc.edu}
\footnote{The author was partially supported by NSF grant DMS-0300356.}
\maketitle
\section{Introduction}
Our aim in this paper is to prove a theorem which implies  a  multiplicative analogue of the Horn conjecture (see ~\cite{ful} for a discussion of the classical case).

To state our theorem we have to
 first recall the setting of two closely related topics: quantum cohomology of the Grassmannians $\Gr(r,n)$ and the multiplicative  eigenvalue problem for $\SU(n)$. In the introduction all schemes are over $\Bbb{C}$.
\subsection{Quantum cohomology and eigenvalue problems}\label{sec1}
Let $I$ be a subset of $\{1,\dots,n\}$ of cardinality $r$. Make the convention that a set $I$ as above is always written in the form $I=\{i_1<\dots<i_r\}$. Let
$$F_{\sssize{\bullet}}:\{0\}=F_0\subset F_1\subset\dots\subset F_n=W$$
be a complete flag in an $n$-dimensional vector space $W$. Define $\Omega_{I}(F_{\sssize{\bullet}})\subseteq\Gr(r,W)$ to be
$$ \{ V \in \Gr(r,W) \mid \rk(V \cap
F_{i_{a}}) \geq a,\ 1\leq a \leq r\}.$$ Denote the cohomology class  of
this subvariety by $\omega_{I}$. The codimension of $\Omega_{I}(F_{\bull})$ is $\codim(\omega_{I})=\sum_{a=1}^r (n-r+a-i_a)$.

Fix a set of  points $\mpp=\{p_1,\dots,p_s\}$ on $\Bbb{P}^1$. Given subsets $I^1,\dots,I^s$ of $\{1,\dots,n\}$ each of cardinality $r$, and a nonnegative integer
$d$, define the Gromov-Witten  number $\langle\omega_{I^1},\dots,\omega_{I^{s}}\rangle_{d}$
 to be, for generic flags $F^j_{\bull}$ on $\cn$ for $j=1,\dots,s$, the number of maps $f:\Bbb{P}^1\to Gr(r,W)$ of degree $d$ such that for each $j=1,\dots,s$,
$f(p_j)\in \Omega_{I^j}(F^j_{\bull})$. If there is an infinite number of such maps, $\langle\omega_{I^1},\dots,\omega_{I^{s}}\rangle_{d}$ is defined to be zero.

Conjugacy classes  in the special unitary group $\SU(n)$ are in one to one correspondence with sequences
of the form $(\delta_1,\dots,\delta_r)$ satisfying
$$\delta_1\geq\dots\geq\delta_n\geq\delta_1-1,\qquad\sum_{b=1}^n \delta_b=0,$$where, to  $(\delta_1,\dots,\delta_r)$, we associate the conjugacy class of the diagonal  matrix with entries $\exp(2\pi i \delta_b)$ for $ b=1,\dots,n$.

For a conjugacy class $\bar{A}$ corresponding to the sequence $(\delta_1,\dots,\delta_r)$,  and a subset $I\subseteq \{1,\dots,n\}$ of cardinality $r$, define
$$\lambda_I(\bar{A})=\sum_{i\in I}\delta_i.$$
The following  theorem  was proven independently by Agnihotri-Woodward ~\cite{AW} and the author ~\cite{b1}. It says that the multiplicative eigenvalue problem for $\SU(n)$ is controlled by quantum Schubert calculus of the Grassmannians $\Gr(r,n)$.
\begin{theorem}\label{oldone}Let $\bar{A}^{(1)},\dots,\bar{A}^{(s)}$ be conjugacy classes in $\SU(n)$. Then,  there exist $A^{(1)},\dots,A^{(s)}$ in $\SU(n)$ with $A^{(j)}$ in the conjugacy class of $\bar{A}^{(j)}$ for $j=1,\dots,s$ and  $A^{(1)}\,  A^{(2)}\, \cdots A^{(s)}\, =\, I$ if and only if:
For any integers $r$, $d$ with $0<r<n$, $d\geq 0$ and subsets
 $I^1,\dots,I^s$ of $\{1,\dots,n\}$ each of cardinality $r$, such that
$\langle\omega_{I^1},\dots,\omega_{I^{s}}\rangle_{d}\neq 0$, the following inequality holds:
\begin{equation}\label{I1}
\sum_{j=1}^{s}\lambda_{I^j}(\bar{A}^{(j)})\leq d.
\end{equation}
\end{theorem}
\subsection{The main results}\label{mt}
We would like to understand (inspired by the classical Horn problem), the condition $\langle\omega_{I^1},\dots,\omega_{I^{s}}\rangle_{d}\neq 0$ from Theorem ~\ref{oldone} in terms of the multiplicative eigenvalue problem for $\SU(r)$. To formulate our result we introduce some notation:

 For $I=\{i_1<\dots<i_r\}\subset\{1,\dots,n\}$, define a conjugacy class
$\beta(I)=(\beta_1,\dots,\beta_r)$ for $\SU(r)$ as follows: First
 define ${\Lambda}(I)=(l_1,\dots,l_r)$ where $l_a=\frac{n-r+a-i_a}{n-r}$ for $a=1,\dots,r$. Let
$c=\frac{1}{r}\sum_{a=1}^r l_a$ and finally, let $\beta_a=l_a-c$ for $a=1,\dots,r$.

The center of $\SU(r)$ acts on the conjugacy classes of elements in $\SU(r)$. To make this explicit, let $\zeta_r=\exp(\frac{2\pi i}{r})\in \Bbb{C}$. To $\zeta_r$ one can associate a generator of the center of $\SU(r)$, namely the diagonal matrix with $\zeta_r$ on the diagonal. Given a conjugacy class $\Delta$ for $\SU(r)$ we have a natural conjugacy class $\zeta_r\Delta$ for $\SU(r)$. Explicitly, if $\Delta=(\delta_1,\dots,\delta_r)$ then
$$\zeta_r\Delta=(\delta_2+\frac{1}{r},\dots,\delta_r+\frac{1}{r},\delta_1+\frac{1}{r}-1).$$
We now state the main theorem of this paper which relates the nonvanishing of
 a Gromov-Witten number to the product of unitary matrices problem (see Corollary ~\ref{viv} for a symmetric form)
:
\begin{theorem}\label{maintheorem} Let $I^1,\dots,I^s$  be subsets of $\{1,\dots,n\}$ each of cardinality $r$ and $d$ a nonnegative integer such that
$$\sum_{j=1}^s \codim(\omega_{I^j})=r(n-r)+dn.$$
The following are equivalent:
\begin{enumerate}
\item $\langle\omega_{I^1},\dots,\omega_{I^{s}}\rangle_{d}\neq 0$.
\item There exist $A^{(1)},\dots,A^{(s)}$ in $\SU(r)$ satisfying
\begin{itemize}
\item $A^{(1)}\,  A^{(2)}\, \cdots A^{(s)}\, =\, I$.
\item $A^{(1)}$ is in the conjugacy class corresponding to $\zeta_r^d{\beta}(I^1)$
 and for $j=2,\dots,s$, $A^{(j)}$ is in the conjugacy class corresponding to ${\beta}(I^j).$
\end{itemize}
\item There exists a (special) unitary local system $\mathcal{L}$ on $\pone-\mpp$ such that the local monodromy of $\mathcal{L}$ at $p_1$ is $\zeta_r^d{\beta}(I^1)$, and  the monodromy at $p_j$ for $j\geq 2$ is ${\beta}(I^j).$
\end{enumerate}
\end{theorem}
The equivalence $(2)\Leftrightarrow(3)$ follows immediately from the description of the fundamental group of $\pone-\mpp$.

The implication $(1)\Rightarrow (3)$ is a remarkable way of producing unitary local systems on $\pone-\mpp$ from non-vanishing Gromov-Witten numbers. This can also be obtained from the work of E. Witten ~\cite{witten} and S. Agnihotri ~\cite{agni} as will be explained in Section ~\ref{agnihotri}. 

Theorem ~\ref{maintheorem} along with Theorem ~\ref{oldone}, can be used to  obtain a generalization (not
conjectured before) of A. Horn's 1962 conjecture ~\cite{Horn} on the eigenvalues of sums of Hermitian matrices
(see Section ~\ref{consequence}). The multiplicative generalization is concerned with an inductive
characterization of the possible eigenvalues of a product of unitary matrices. The original 1962 conjecture of
Horn was proved by the combined works of A. Klyachko, A. Knutson and T. Tao ~\cite{Klyachko}, ~\cite{KT} (see
~\cite{ful} for a history of this problem). 

We prove a more general theorem than Theorem ~\ref{maintheorem} (see Theorem ~\ref{MainTe}) which gives information (Corollary ~\ref{info}, (1)$\Leftrightarrow$(2)) on the smallest power of $q$ in a
quantum product of Schubert varieties in Grassmannians with an arbitrary number of factors. The case of two factors for an arbitrary $G/P$ was  considered by W. Fulton and C. Woodward in ~\cite{fulW}.

The quantum analogue of the saturation theorem of Knutson and Tao ~\cite{KT} is stated in a non geometric form in Section ~\ref{qsat}. The geometric form of this theorem
will be given in ~\cite{bQ2}.

Our methods give transversality statements in quantum Schubert calculus in any characteristic (see Section ~\ref{gamma}). The interesting problem of the maximum possible number of real or $p$-adic solutions to a ``quantum Schubert enumerative problem'' remains open (but see ~\cite{sottile3}, ~\cite{vakil}).

\subsection{An overview of the methods}
It is standard that one  can view $\langle\omega_{I^1},\dots,\omega_{I^{s}}\rangle_{d}$ form Section ~\ref{sec1} also as the number of subbundles (if finite and zero otherwise) $\mv$ of $\mw=\mathcal{O}^n$ so that for each $p\in\mpp$, $\mv_p\in \Omega_{I^j}(F^j_{\bull})$. To get an inductive grip on this situation, we would like to replace $\mw$ by $\mv$. However $\mv$ is of degree $-d$ which is not necessarily $0$. It can however be shown to be evenly split (see Section ~\ref{rembrandt} and Lemma ~\ref{openn}). This motivates us to carry out a generalization of Gromov-Witten numbers (Section ~\ref{temina}). 

With this inductive framework in place, the strategy for the  proofs is very similar to those in ~\cite{jag} (where many of the arguments in this paper appear in a simpler situation): The tangent space technique is modified so that it applies in the context of space of maps of $\pone$ to a homogenous space. We use standard properties of Quot schemes to do the tangent space calculations. We also use the  general position techniques from ~\cite{jag} (see ~\cite{jag}, Introduction).

 There are some additional difficulties  in the quantum situation arising from the nature of maps between vector bundles on $\pone$ (eg. image of a morphism of vector bundles may not be a subbundle). We also use techniques inspired by the theory of parabolic bundles to overcome these difficulties. 
\subsection{Conventions}\label{conventions}
All schemes in this paper are assumed to be  finite type over an algebraically closed base field $\kappa$ of arbitrary characteristic.
\begin{enumerate}
\item Fix a  finite collection of points $\mpp=\{p_1,\dots,p_s\}$ on $\Bbb{P}^1$.
\item  A vector bundle $\mv$ on $\pone \times\ X$ is said to be a $(d,r)$-bundle if for each $x\in X$, $\mv_x$ is a vector bundle on $\pone$ of degree $-d$ and rank $r$.
\item A morphism $\mv\ \to\ \mw$ of locally free sheaves on a scheme $X$ is said to have rank $r$ if the cokernel is a locally free sheaf of rank $\rk(\mw)-r$.
\item If $\mv$ is a vector bundle on a scheme $X$, the contravariant functor schemes/$X$ to (sets) given $T\leadsto$ the set of  complete filtrations by subbundles $F_{\bull}$:
$$0\subsetneq F_1\subsetneq F_2\subsetneq\dots\subsetneq F_r=\mv_T$$
of $\mv_T$ ($\mv_T$ denotes the pullback of $\mv$ to $T$), is representable  by a flag variety $\Flag(\mv)$ which is smooth over $X$.
\item Denote the set $\{1,\dots,r\}$ by $[r]$.
\end{enumerate}
\subsection{Acknowledgements}
I thank A. Buch, W. Fulton, F. Sottile, C. Woodward and A. Yong for useful
communication. Woodward pointed out the related work of Witten ~\cite{witten} and Agnihotri ~\cite{agni}
and that it should give a different proof of $(1)\Rightarrow (2)$ in Theorem ~\ref{maintheorem} (see Section ~\ref{agnihotri}). I thank Xiaowei Wang for giving me a copy of Agnihotri's 1995 Ph.D thesis ~\cite{agni}.

\section{Formulation of the Main result}
\subsection{Evenly split bundles on $\pone$}\label{rembrandt}
 A $(D,n)$-vector bundle (see Conventions ~\ref{conventions}) $\mw$ on $\pone$ is said to be evenly split (ES)
if $\mathcal{W}=\oplus_{i=1}^n\mathcal{O}_{\Bbb{P}^1}(a_i)$  with $|a_i-a_j| \leq 1$ for $0<i<j\leq n$. It is easy to see that $\mw$ is ES if and only if  $\Hone(\pone,\send(\mw))=0$.

Let $D$, and $n$ be integers with $n>0$. It is easy to show that upto  isomorphism, there is a unique  ES -bundle of degree $-D$ and rank $n$ on $\pone$. We denote this bundle by $\mz_{D,n}$.

 Let $\mw$ be a bundle on $\pone$. Define {$\mfm(d,r,\mw)$} to be the moduli space of $(d,r)$-subbundles of  $\mw$. This can be obtained as an open subset of the Quot scheme of quotients of $\mw$ of degree $d-D$ and rank $n-r$. In the notation of ~\cite{potier}, $\mfm(d,r,\mw)$ is the open subset of  $\operatorname{Hilb}^{n-r,d-D}(\mw)$ formed by points where the quotient is locally free.

If $D,d$ are integers and $0\leq r\leq n$, define {$\mfm(d,r,D,n)$}$=\mfm(d,r,\mz_{D,n})$.
\begin{defi}
For $r$, $m$ positive integers and $d$, $b\in\Bbb{Z}$ define
$$\cee(d,r,b,m)=\chi(\pone,\shom(\mz_{d,r},\mz_{b,m }))=rm \ +\ dm -br.$$
\end{defi}
\begin{proposition}\label{bigone}
$\mfm(d,r,D,n)$ is smooth and irreducible of dimension $\cee(d,r,D-d,n-r)$. The subset of $\mfm(d,r,D,n)$ formed by ES-subbundles $\mv\subseteq\mz_{D,n}$ such that $\mz_{D,n}/\mv$ is also ES, is open and dense in $\mfm(d,r,D,n)$.
\end{proposition}
The proof of  Proposition ~\ref{bigone} will be given in  Section ~\ref{proofbigone}.
\subsection{Complete Flags}\label{old}
For a bundle $\mw$ on $\pone$, define
$$\Fl(\mw)=\prod_{p\in \mpp}\ \Flag(\mw_p).$$
If $\me\in \Fl(\mw)$, we will assume that it is written in the form
$\me=\prod_{p\in\mpp}\ E^p_{\bull}$.
More generally if $X$ is a scheme and $\mw$ is a vector bundle on $\pone\times\  X$, let $\Fl(\mw)$ be the scheme over $X$, whose fiber over $x\in X$ is $\Fl(\mw_x)$ as defined before.

For a bundle $\mw$ on $\pone$, a subbundle $\mv\subseteq \mw$ and a collection of flags $\me=\prod_p E^p_{\bull}\in \Fl(\mw)$ we have associated induced complete flags on $\mv$ and on $\mq=\mw/\mv$ at points of $\mpp$. We denote these by $\me(\mv)=\prod_p E^p_{\bull}(\mv)\in \Fl(\mv)$ and $\me(\mq)=\prod_p E^p_{\bull}(\mq)\in \Fl(\mq)$.

\subsection{Schubert states and Generalized Gromov-Witten numbers}\label{temina}
A Schubert state is a $5$-tuple $\mi=(d,r,D,n,I)$ where
$d$, $D$, $r$ and $n$ are integers, $n\geq r\geq 0$ and $I$ is an assignment to each $p\in\mpp$ a subset $I^p$ of $[n]=\{1,\dots,n\}$ of cardinality $r$.  We will use the notation
$I^{p}=\{i^p_1<\dots<i^p_r\}$ for $p\in \mpp$.

Let $\mw=\mz_{D,n}$, $\me\in \Fl(\mw)$ a generic point and $\mi=(d,r,D,n,I)$ a Schubert state. For $p\in\mpp$ let $\pi_p: \mfm(d,r,\mw)\to \Gr(r,\mw_p)$ be the natural map (the fiber of the subbundle at $p$). Define $\braI$ to be the number of points  in the intersection (if finite and $0$ otherwise)
\begin{equation}\label{omegai}
\bigcap_{p\in\mpp}\pi_p^{-1}[\Omega^o_{I^{p}}(E^p_{\bull})]\subseteq\mfm(d,r,\mw).
\end{equation}
For a Schubert state $\mi$ as above, define
\begin{equation}\label{expected}
\dimii=\dim\mfm(d,r,D,n)-\sum_{p\in\mpp}\codim(\omega_{I^{p}}).
\end{equation}
\begin{remark} If $D=0$ this corresponds to the usual definition of Gromov-Witten numbers. The generalization above  makes induction arguments possible (to pass from the pair $(D,n)$ to the pair $(d,r)$). On the other hand, these new numbers can be recovered from the ``usual'' Gromov-Witten numbers (Corollary ~\ref{vinci}).
\end{remark}
For $p\in\mpp$, the group $\Gl(\mw_p)$ acts transitively on  $\Gr(r,\mw_p)$. Hence by a theorem of Kleiman (see ~\cite{int}, \S B.9.2), the dimension of Intersection ~\ref{omegai} is given by $\dimii$ (with $\mf$ generic).

It is important for us to include study the cases when the Intersection ~\ref{omegai} is nonempty. We will say that {\it $\mi$ is not  null} if the intersection ~\ref{omegai} is nonempty (possibly infinite) for generic $\me\in \Fl(\mw)$.

If the base field $\kappa$ is of characteristic $0$, it follows from Kleiman's transversality theorem that for any Schubert state  $\mi=(d,r,D,n,I)$ such that the expected dimension given by Equation ~\ref{expected}  (viz. $\dimii$) is zero and generic $\me\in\Fl(\mw)$,  Intersection ~\ref{omegai} is a reduced zero dimensional scheme. In Section  ~\ref{gamma}, we will show that this property holds for any algebraically closed field.

\subsection{Shift operations}\label{termin}
The following two results on shift operations are proved in Section ~\ref{PP}.
\begin{lemma}\label{jkl}
Let $\mi=(d,r,D,n,I)$ be a Schubert state. Let $\mj=(d+r,r,D+n,n,I)$. Then,
\begin{enumerate}
\item The expected dimensions of the intersections corresponding to  $\mi$ and $\mj$ are the same. That is, $\dimii=\dim\mj$.
\item  $\mi$ is not null $\Leftrightarrow$ $\mj$ is not null.
\item $\langle \mi\rangle=\langle \mj\rangle.$
\end{enumerate}
\end{lemma}
 Let $\mi=(d,r,D,n,I)$ be a Schubert state and $p\in \mpp$. Define
a new shifted Schubert state $\sS(p)(\mi)=(\tilde{d},r,D-1,n,J)$ where $\tilde{d}$ and $J$ are given  as follows: $J^q=I^q$ if
$q\in \mpp-\{p\}$, and letting $I^{p}=\{i_1<\dots<i_r\},$
\begin{enumerate}
\item If $i_1>1$, let $J^{p}=\{i_1-1<\dots <i_r-1\}$ and $\tilde{d}=d$.
\item If $i_1=1$, let $J^{p}=\{i_2-1<\dots <i_r-1<n\}$ and $\tilde{d}=d-1$.
\end{enumerate}
The following proposition  is related to the ``action'' of $n$th roots of unity on the Quantum cohomology of the
Grassmannian $\Gr(r,n)$ ~\cite{AW}. The geometry of this action of the center of $\SL(n)$ appeared in ~\cite{b2}.

\begin{proposition}\label{shiftt1}
With notation as above,
\begin{enumerate}
\item The expected dimensions of the intersections corresponding to  $\mi$ and $\mj$ are the same. That is, $\dimii=\dim\mj$
\item  $\mi$ is not null $\Leftrightarrow$  $\sS(p)\mi$ is not null.
\item $\langle\mi\rangle=\langle \sS(p)(\mi)\rangle.$
\end{enumerate}
\end{proposition}
The above results have the following useful corollary:
\begin{corollary}\label{vinci}
Let $\mi=(d,r,D,n,I)$ be a Schubert state. Then one can determine a Schubert state $\mj=(\td,r,n,0,J)$ such that
\begin{enumerate}
\item
$\langle\mi\rangle=\langle\mj\rangle$. The right hand side is an ``usual'' Gromov-Witten number.
\item $\mi$ is not null $\Leftrightarrow$  $\mj$ is not null.
\end{enumerate}
\end{corollary}
\begin{proof}
Using Lemma ~\ref{jkl}, assume $0\leq D<n$. Now make a $D$-fold application of the shift operator $\sS(p)$.
\end{proof}
\subsection{Statement of the main result}
\begin{theorem}\label{MainTe}
Let $\mi=(d,r,D,n,I)$ be a Schubert state such that
$$\dimii\ =\ \dim\mfm(d,r,D,n)-\sum_{p\in\mpp}\codim(\omega_{I^{p}})\geq 0$$
The following are equivalent:
\begin{enumerate}
\item[(A)] $\mi$ is not null.
\item[(B)] Given a non null Schubert state of the form $\mk=(\td,\tr,d,r,K)$ with
$0<\tr<r$, the following inequality holds:
\begin{equation}\nonumber\tag{$\dagger_{\mk}^{\mi}$}
-\cee(\td,\tr,D-d,n-r)\ +\ \sum_{p\in\mpp}\sum_{a\in K^{p}}(n-r+ a-i^p_a)\leq 0.
\end{equation}
Explicitly, Inequality $(\dagger_{\mk}^{\mi})$ is the following,
$$-\td(n-r)\ +\  \tr (D-d)-\tr (n-r)\ +\ \sum_p\sum_{a\in K^{p}}(n-r+ a-i^p_a)\leq 0$$
\item [(C)]Given a Schubert state of the form $\mk=(\td,\tr,d,r,K)$ with
$0<\tr<r$ and $\langle \mk\rangle\neq 0$, Inequality ($\dagger_{\mk}^{\mi}$)
holds.
\item [(D)]Given a Schubert state of the form $\mk=(\td,\tr,d,r,K)$ with
$0<\tr<r$ and $\langle \mk\rangle=1$, Inequality ($\dagger_{\mk}^{\mi}$)
holds.
\end{enumerate}
\end{theorem}
\begin{remark}
The proofs of $(A)\Leftrightarrow(B)\Leftrightarrow(C)$ are independent of transversality statements (Section ~\ref{gamma}), and do not invoke the results of Section ~\ref{PP}. The transversality result of Section ~\ref{gamma} is deduced from the equivalence of $(A)$, $(B)$ and $(C)$.

The transversality statement in Section ~\ref{gamma} is used to show the equivalence of $(D)$ with the other three conditions. 
\end{remark}
\section{First Applications of Theorem ~\ref{MainTe}}
\subsection{Proof of Theorem ~\ref{maintheorem}}\label{thissection}\label{loungeback}
Let $\kappa=\Bbb{C}$. We introduce some notation for this section: Conjugacy classes in $\SU(n)$ are parameterized by points in the $n-1$ dimensional simplex:
$$\Delta(n)=\{(\delta_1,\dots,\delta_n)\in \Bbb{R}^n\mid \delta_1\geq\dots\geq\delta_n\geq\delta_1-1\ ,\sum_{b=1}^n \delta_b=0\}.$$
To an  element $(\delta_1,\dots,\delta_n)\in \Delta(n)$, we associate the conjugacy class of the diagonal  matrix with entries $\exp(2\pi i \delta_b)$,  $b=1,\dots,n$ on the diagonal.

Now define  $\Gamma(n,s)\subseteq \Delta(n)^s$ be the set of $(\Delta^1,\dots,\Delta^s)\in \Delta(n)^s$ such that there exist $A^{(j)}\in \SU(n)$ in the conjugacy class $\Delta^j$ for $j=1,\dots,s$ satisfying the equality $A^{(1)} A^{(2)}\dots A^{(s)}=I$. Theorem ~\ref{oldone} gives a description of $\Gamma(n,s)$ in terms of inequalities.
We note the following corollaries of Theorem ~\ref{MainTe}.
\begin{corollary}\label{initial}
Let $\mi=(d,r,D,n,I)$ be a Schubert state such that $d=0$, and
$\dimii=0.$
Then the following are equivalent
\begin{enumerate}
\item  $\langle\mi\rangle \neq 0$.
\item $({\beta}(I^{p_1}),{\beta}(I^{p_2}),\dots,{\beta}(I^{p_s}))\in \Gamma(r,s)$.
\end{enumerate}
\end{corollary}
\begin{proof}
Let ${\beta}(I^{p_j})=({\nation}^j_1,\dots,{\nation}^j_r)$ for
$j=1,\dots,s$. Let $b^j_a=\frac{n-r+a-i^j_a}{n-r}$, $a=1,\dots,r$. From the definitions, we have ${\nation}^j_a=b^j_a-c_j$ where $rc_j=\sum_{a=1}^r b^j_a$.

From Theorem ~\ref{MainTe} $(A)\Leftrightarrow(C)$, the non nullness of $\mi$ is equivalent to a set of inequalities indexed by Schubert states of the form $\mk=(\td ,\tr ,0,r,K)$ satisfying $\langle\mk\rangle \neq 0$. Similarly, according to Theorem ~\ref{oldone}, (2) holds iff a system of inequalities indexed by the same set of $\mk$ holds. We just have to check the inequalities are the same. The inequality corresponding to Theorem ~\ref{MainTe}, $(\dagger_{\mk}^{\mi}$) is:
$$\frac{1}{\tr}(-\td +\sum_{j=1}^s\sum_{a\in K^{p_j}}b^j_a)\ +\  \frac{D}{n-r}-1\leq 0 $$
which is the same as
\begin{equation}\label{feirst}
\frac{1}{\tr}(-\td +\sum_{j=1}^s\sum_{a\in K^{p_j}}{\nation}^j_a)\ +\  \sum_j c_j +\frac{D}{n-r}-1\leq 0 .
\end{equation}
Using the hypothesis $\dimii=0$, we have
$$r(n-r)[\sum_j c_j \ +\ \frac{D}{n-r}-1]= \sum_{p\in \mpp} \codim(\omega_{I^{p}}) \ +\  Dr-r(n-r) =0.$$
Therefore Inequality ~\ref{feirst} is
$$\frac{1}{\tr}(-\td \ +\ \sum_{j=1}^s\sum_{a\in K^{p_j}}{\nation}^j_a)\leq 0 .$$
The above inequality  is the same as the one corresponding to Theorem ~\ref{oldone} for $\mk$ and the corollary is proved.
\end{proof}
Consider the operator $T(r,n)$ acting on subsets $I=\{i_1<\dots <i_r\}$ of
$[n]$ which takes $I$ to
\begin{enumerate}
\item $\{i_1-1<i_2-1<\dots<i_r-1\}$ if $i_1>1$.
\item $\{i_2-1<\dots< i_r-1<n\}$ if $i_1=1$.
\end{enumerate}
Informally, $T(r,n)I=I-1$ with $0$'s replaced by $n$. It is immediate that  ${\beta}(T(r,n)I)={\beta}(I)$ if $i_1\neq 1$ and $\zeta_r{\beta}(I)$ if $i_1=1$. The following {\bf includes  Theorem ~\ref{maintheorem} as a special case:}
\begin{corollary}\label{terminator}
Let $\mi=(d,r,D,n,I)$ be a Schubert state such that
$\dimii=0$.
Then,  $\langle\mi\rangle \neq 0$ if and only if
$ \Leftrightarrow (\zeta_r^d{\beta}(I^{p_1}),{\beta}(I^{p_2}),\dots,{\beta}(I^{p_s}))\in \Gamma(r,s)$.
\end{corollary}
\begin{proof}
We have the following two equivalences:
$\langle\mi\rangle\neq 0$ iff  $\langle(d+r,r,D+n,I)\rangle\neq 0$ and
$(\zeta_r^d{\beta}(I^{p_1}),{\beta}(I^{p_2}),\dots,{\beta}(I^{p_s}))\in \Gamma(r,s)$ iff
$(\zeta_r^{d+r}{\beta}(I^{p_1}),{\beta}(I^{p_2}),\dots,{\beta}(I^{p_s}))\in \Gamma(r,s)$.
The second equivalence is obvious and the first one is from Lemma ~\ref{jkl}.
We can therefore assume that $0\leq d< r$ without loss of generality.

Let $I^{p_1}=\{i_1<\dots<i_r\}$, we consider $\mj=\sS(p_1)^{i_d}\mi=(0,r,D-i_d,n,J)$ (see Section ~\ref{termin} for definition of $\sS(p)$ for $p\in \mpp$). Clearly $J^{p_1}=T(r,n)^{i_d}I^{p_1}$ and $J^{q}=I^{q}$ for $q\neq p$, $q\in\mpp$.

From  Proposition ~\ref{shiftt1}, we have $\langle\mi\rangle\neq 0$ if and only if
$\langle \mj\rangle \neq 0$. Also, ${\beta}(J^{p_1})=\zeta^d{\beta}(I^{p_1})$ and ${\beta}(J^{p_j})={\beta}(I^{p_j})$ for $j=2,\dots,s$. We can therefore appeal to Corollary ~\ref{initial} (for the Schubert state $\mj$) and the Corollary is proved.
\end{proof}
For any integers $(a_1,\dots,a_s)$ such that $\sum_{j=1}^s a_j$ is divisible by $r$, the map $\Delta(r)^s\to \Delta(r)^s$ which takes
$(\Delta^1,\dots,\Delta^s)$ to $(\zeta_r^{a_1}\Delta^1,\dots,\zeta_r^{a_s}\Delta^s)$ clearly preserves the set $\Gamma(r,s)$. We therefore have the following  symmetric form of Theorem ~\ref{maintheorem}:
\begin{corollary}\label{viv}
Let $I^1,\dots,I^s$ be subsets of $[n]$ of cardinality $r$ each and $d\geq 0$ a positive integer such that $\sum_{j=1}^s \codim(\omega_{I^j})=r(n-r)\ +\ dn.$ Also assume that we are given $(a_1,\dots,a_s)\in \Bbb{Z}^s$ such that $\sum_{j=1}^s a_j-d$ is divisible by $r$. Then,
$\langle\omega_{I^1},\dots,\omega_{I^{s}}\rangle_{d}\neq 0$ if and only if $(\zeta_r^{a_1}{\beta}(I^1),\zeta_r^{a_2}{\beta}(I^2),\dots,\zeta_r^{a_s}{\beta}(I^s))\in \Gamma(r,s).$
\end{corollary}
\subsection{Consequences for eigenvalue problems}\label{consequence}
Theorem ~\ref{oldone} together with Theorem ~\ref{maintheorem} (Corollary ~\ref{terminator}) implies an inductive
characterization of the  list of inequalities for the multiplicative eigenvalue problem. This characterization
does not involve quantum cohomology, and can be considered to the multiplicative analogue of Horn's original
conjecture.

For integers $0<r<n$, define  $A(r,n,s)$ to be the set of tuples $(d,r,n,I^1,\dots,I^s)$ where $d$ is a nonnegative  integer, $I^j\subseteq [n]$ for $j=1,\dots,n$ are subsets of cardinality $r$ each  satisfying the condition
$$\sum_{j=1}^s\sum_{a=1}^r (n-r+a-i^j_a)=r(n-r)\ +\ dn.$$

We inductively define subsets $\tilde{\Gamma}(n,s)\subseteq \Delta(n)^s$ and $B(r,n,s)\subseteq A(r,n,s)$ where $r,n$ are integers such that $0<r<n$ as follows \begin{itemize}
\item $\tilde{\Gamma}(1,s)=\Delta(1)^s$ and $B(1,n,s)=A(1,n,s)$ for all positive integers $n$.
\item Let $n$ be a positive integer. Assume $B(r,k,s)$ has been defined whenever $r<n$ and $\tilde{\Gamma}(k,s)$ has been defined whenever $k<n$.
Define $\tilde{\Gamma}(n,s)$ and for $m>n$, the set  $B(n,m,s)$  successively as follows:
\begin{enumerate}
\item $(\Delta^1,\dots,\Delta^s)\in \tilde{\Gamma}(n,s)$ where $\Delta^j=(\delta^j_1,\dots,\delta^j_n)$ for $j=1,\dots,s$ if and only if: For any integers $\tr $, $d$ with $n> \tr > 0$ and subsets
$I^1,\dots,I^s$ of $[n]$ each of cardinality $\tr $, such that
$(d,\tr ,n,I^1,\dots,I^s)\in B(\tr ,n,s)$, the following inequality holds:
$$\sum_{j=1}^{s}\sum_{b\in I^j}\delta^j_b\leq d.$$
\item $(d,n,m,I^1,\dots,I^s)\in A(n,m,s)$ is an element of $B(n,m,s)$
 if and only if
$$(\zeta_n^d{\beta}(I^1),{\beta}(I^2),\dots,{\beta}(I^s))\in \tilde{\Gamma}(n,s).$$
\end{enumerate}
\end{itemize}
The following corollary is immediate from Corollary ~\ref{terminator} and  Theorem ~\ref{oldone}.
\begin{corollary}
\begin{enumerate}
\item For any positive integer $n$, $\Gamma(n,s)=\tilde{\Gamma}(n,s)$.
\item Suppose we are given  $0<r<n$ a positive integer, an integer $d\geq 0$
 and subsets $I^1,\dots,I^s$ of $[n]$ of cardinality $r$ each. Then,
$\langle\omega_{I^1},\dots,\omega_{I^{s}}\rangle_{d}\neq 0$ iff
 $(d,r,n,I^1,\dots,I^s)\in B(r,n,s)$.
\end{enumerate}
\end{corollary}
\subsection{Lowest powers of $q$}\label{lounge}
We want to give a criterion, expressed in terms of small quantum cohomology, for a Schubert state of the form $\mi=(d,r,n,0,I)$ to be non null.
\begin{lemma}\label{qcr}
Given $\omega_J$, the cohomology class of a Schubert variety in $\Gr(r,n)$ and an integer $d> 0$, there exist $K_1,\dots, K_{\ell}$ subsets of $[n]$ each of cardinality $r$ such that in the (small quantum) product
$$\omega_J\star\omega_{K_1}\star\dots\star\omega_{K_{\ell}},$$
 the coefficient of $q^d\omega_L$ is nonzero for some $L$.
\end{lemma}
\begin{proof} We easily see that we need to take care of only the case $d=1$.
Multiply by the Poincare dual of $\omega_J$ to reduce to the case of $\omega_J = $ class of a point. So we need to exhibit a $K$ so that $\omega_K\star[pt]$
has a term in which  $q$ appears to the degree 1. It is easy to see that the codimension $1$ Schubert variety does this job (use Bertram's Pieri formula ~\cite{bert}).
\end{proof}
\begin{proposition}\label{lowestq} If $\mi=(d,r,D,n,I)$ and $D=0$, then  $\mi$ is not null if and only if  in the (small) quantum product $$\omega_{I^{p_1}}\star\omega_{I^{p_2}}\star\dots\star\omega_{I^{p_s}},$$
the coefficient of $q^c\omega_J$ is nonzero for some $c$ satisfying $0\leq c\leq d$, and $J$ a subset of $[n]$ of cardinality $r$.
\end{proposition}
\begin{proof}
Assume that the quantum star product has a term $q^c\omega_J$ with $c\leq d$. According to Lemma ~\ref{qcr} we can find $K_1,\dots,K_m$ subsets of
$[n]$ each of cardinality $r$ such that
$$\langle\omega_{I^{p_1}},\dots,\omega_{I^{p_s}},\omega_{K_1},\dots,\omega_{K_m}\rangle_d\neq 0.$$
It is now immediate that $\mi$ is not null.

To go the other way, we let $\mw$ be an ES $(D,n)$-bundle on $\pone$, $\mf\in\Fl(\mw)$ such that $\Omega^o(\mi,\mw,\mf)$ is of the expected dimension. Let $m=\dimii$. Pick a large collection of points  $\mq=\{q_1,\dots,q_{\ell}\}$ such that $$\Gamma:\mfm(d,r,\mw)\to \prod_{q\in\mq}\Gr(r,\mw_q)$$ is an embedding (finite
 will do). By intersecting the image of $\Gamma$ with $(\sum_{a=1}^{\ell} D_a)^m$ where $D_a$ is the codimension $1$ Schubert variety on  $\Gr(r,\mw_{q_a})$ we find that
 $$\langle\omega_{I^{p_1}},\dots,\omega_{I^{p_s}},\omega_{K_1},\dots,\omega_{K_{\ell}}\rangle_d\neq 0$$
for some choice of $\omega_{K_i}$ for $i=1,\dots,{\ell}$.

Using the associativity of the quantum product, one finds that in
$$\omega_{I^{p_1}}\star\omega_{I^{p_2}}\star\dots\star\omega_{I^{p_s}},$$
the coefficient of $q^c\omega_J$ is nonzero for a $0\leq c\leq d$ and $J$ a subset of $[n]$ of cardinality $r$.
\end{proof}

Recall the operators $T(r,n)$ and $\sS(p)$ for $p\in\mpp$ defined in Section ~\ref{thissection} and  Section ~\ref{termin} respectively. We then have the following:
\begin{corollary}\label{info}
Let $I^1, \dots I^s$ be subsets of $[n]$, each of cardinality $r$.
Let $d\geq 0$ be an integer. Write $d=qr+b$ with $0\leq b<r$ and $(q,b)\in\Bbb{Z}^2$. let $\tilde{I}^1=T(r,n)^{i^1_b} I^1=I^1-i^1_b$ (with $0$'s replaced by $n$'s and $i^1_0$ defined to be $0$).
The following are equivalent
\begin{enumerate}
\item[(1)]  In $\omega_{I^1}\star\dots\star\omega_{I^s},$ a term of the form $q^c\omega_J$ with  $0\leq c\leq d$ appears with a non-zero coefficient.
\item[(2)] $\sum_{j=1}^s\codim(\omega_{I^j})\leq dn\ +\  r(n-r)$
 and for every choice of $(\td,\tr,K^1,\dots,K^s)$ where $\td\geq 0$ and  $r>\tr>0$ are integers and $K^1,\dots,K^s$ are subsets of  $[r]=\{1,\dots,r\}$ each of cardinality $\tr$ such that $\langle\omega_{K^1},\dots,\omega_{K^s}\rangle_{\td}=1,$ the inequality
$$\sum_{a\in K^1}(n-r+ a-\tilde{i}^1_a)+ \sum_{j=2}^{s}\sum_{a\in K^j}(n-r + a-i^j_a)\leq \chi(\td,\tr,-(qn+i^1_b),n-r)\ $$
is valid where
$ \chi(\td,\tr,-(qn+i^1_b),n-r)\   =\ \td (n-r)+\tr (qn+i^1_b)+\tr (n-r).$

\end{enumerate}
\end{corollary}
\begin{proof}Let $I(p_j)=I^j$ and $\mi=(d,r,0,n,I)$. By Proposition ~\ref{lowestq}, the condition in (1) is that $\mi$ is non null.
Let $\mj=\sS(p_1)^{qn+i_b}\mi=(0,r,-(qn+i_b),n,J)$ for some $J$.
By Lemma ~\ref{jkl} and Proposition ~\ref{shiftt1}, $\dimii=\dim\mj$ and $\mi$ is non null if and only if $\mj$ is non null. By Theorem ~\ref{MainTe}, $\mj$ is non null if and only if the conditions in $(2)$ hold (note that $J(p_1)=\tilde{I}^1$ and $J(p_j)=I^j$ for $j=2,\dots,s$).
\end{proof}

\section{Relation to work of Witten and Agnihotri}\label{agnihotri}

It should be pointed out that the statement of $(1)\Rightarrow (2)$ in Theorem ~\ref{maintheorem} has not appeared before. Computations to guess the form of these conditions were made by Buch and Fulton (unpublished). As pointed out to us by Woodward, one may also obtain this implication  from the results of  Witten ~\cite{witten}, which were proven mathematically by Agnihotri in his thesis ~\cite{agni}. We explain this in this section.

Let $b\in\pone-\ \mpp$. Represent $\pi_1(\pone-\ \mpp,\ b)$ in the standard manner: Let $\gamma_j$ be the loop based at $b$ that loops around $p_j$ for $j=1,\dots,s$. Then, the fundamental group $\pi_1(\pone-\ \mpp,\ b)$ is the free group on $\gamma_1,\gamma_2,\dots,\gamma_s$ modulo the smallest normal subgroup generated by the word $\gamma_1\cdot\gamma_2\cdot\dots\cdot\gamma_s$.

Let $G$ be a group. A representation $\rho:\pi_1(\pone-\ \mpp,\ b)\to G$  is the same as choice of elements $A_j\in G$ for $j=1,\dots,s$ with $A_1A_2\dots A_s=I$ (by the association $A_j\leadsto \rho(\gamma_j)$).

Use notation from Theorem ~\ref{maintheorem} and assume that $\sum_{j=1}^s \codim(\omega_{I^j})=r({\emma})+dn$. The theorem of Witten (as in Agnihotri's  thesis) gives an expression
$$\langle\omega_{I^1},\dots,\omega_{I^{s}}\rangle_{d}= h^0(\mathcal{M}(I^1,\dots,I^s,d),\Theta)$$
where $\mathcal{M}=\mathcal{M}(I^1,\dots,I^s,d)$ is the moduli space of semistable parabolic bundles on $\pone$
whose underlying bundle has degree $-d$, rank $r$ and parabolic structure at $p_1,\dots,p_s$ (in the notation of
Mehta and Seshadri's paper ~\cite{MS}). The parabolic weight at $p_j$ given by ${\Lambda}(I^j)$ (defined in Section
~\ref{mt}), and $\Theta$ is an ample line bundle on $\mathcal{M}$.

One therefore deduces that if $\langle\omega_{I^1},\dots,\omega_{I^{s}}\rangle_{d}\neq 0$, then $\mathcal{M}\neq \emptyset$. However if $\mv\in\mathcal{M}$, the parabolic slope of $\mv$ is $1+\frac{d}{{\emma}}$ which is not zero, so the theorem of Mehta-Seshadri cannot be immediately applied
 to get a unitary representation (with local monodromies given by the parabolic weights ${\Lambda}(I^j)$) of $\pi_1(\pone-\ \mpp,\ b)\to \Unit(r)$. We will deal with this difficulty in an ad hoc manner here.

The addition of a parabolic point $q\in\pone-\ \mpp$ with central weight $c\in [0,1)$ (that is, the sequence of weights at $q$ is $(c,c,\dots,c)\in \Bbb{R}^r$) increases the parabolic slope of $\mv$ by $c$. Twisting $\mv$ by $\mathcal{O}(1)$ increases the parabolic weight by $1$. Therefore write an equation
$$c+a+1+\frac{d}{{\emma}}=0$$
where $a\in\Bbb{Z}$, $c\in [0,1)$. Consider a new parabolic bundle
 where we have added a new parabolic point $q$ with central weight $c$ and twisted $\mv$ by $\mathcal{O}(a)$. These operations do not change semistability. The parabolic slope has now become $0$.

Recall the notation from Section ~\ref{mt}. Let ${\beta}(I^j)=(\beta^j_1,\dots,\beta^j_r)$,  ${\Lambda}(I^j)=(l^j_1,\dots,l^j_r)$ and $c^j=\frac{\sum_{a=1}^r l^j_a}{r}$. From the above reasoning we get matrices $A^{(j)}\in \Unit(r)$ for $j=1,\dots,s$ and a scalar matrix $B$ such that:
\begin{enumerate}
\item For $j=1,\dots,s$, $A^{(j)}$ is conjugate to the diagonal matrix with entries $\exp(2\pi i l^j_a)$ along the main diagonal.
\item $B$ is the central diagonal matrix with entries $\exp(2\pi i c)$.
\item $A^{(1)}\cdot A^{(2)}\cdot\dots\cdot A^{(s)}\cdot B=I$.
\end{enumerate}
We calculate
$$({\emma})r\sum_{j=1}^s c^j\ =\ \sum_{j=1}^s\sum_{a=1}^r l^j_a\ =\ \sum_{j=1}^s \codim(\omega_{I(p_j)})\ =\  dn+r({\emma}).$$
Therefore,
$$\sum_{j=1}^s c_j\ =\ \frac{d({\emma})+dr+r({\emma})}{r({\emma})}\ =\ \frac{d}{r}+\frac{d}{{\emma}}+1\ =\ -(c+a)+\frac{d}{r}.$$
Therefore multiplying matrices $A^{(j)}$ by $\exp(-2\pi i c^j)$ we obtain a relation $$\tilde{A}^{(1)}\cdot \tilde{A}^{(2)}\cdot\dots\cdot \tilde{A}^{(s)}\cdot \tilde{B}=I$$
where $\tilde{A}^{(j)}$ are in $\SU(r)$ and in the conjugacy class of ${\beta}(I^j)$
for $j=1,\dots,s $ and $B$ is a diagonal matrix which corresponds to multiplication by $\exp(2\pi i \tilde{c})$ where
$$\tilde{c}\ =\ c+(-(c+a)+\frac{d}{r})\ =\ -a+\frac{d}{r}.$$
Hence $B=\zeta_r^d$ and we are done because we can multiply $\tilde{A}^{(1)}$ by $B$ and write an equation $\tilde{A}^{(1)}\dots\tilde{A}^{(s)}=I$ where
\begin{enumerate}
\item For $j=2,\dots,s$, $\tilde{A}^{(j)}$ is in the conjugacy class of ${\beta}(I^j)$.
\item $\tilde{A}^{(1)}$ is in the conjugacy class of $\zeta_r^d{\beta}(I^1)$.
\item $\tilde{A}^{(1)}\cdot \tilde{A}^{(2)}\cdot\dots\cdot \tilde{A}^{(s)}\ =\ I$.
\end{enumerate}
\subsection{Quantum saturation}\label{qsat}
The classical part of the story had one other aspect - the saturation theorem. We now state the quantum
generalization of this result. We saw above that in the case $\sum_{j=1}^s \codim(\omega_{I^j})=r({\emma})+dn$,
$$\langle\omega_{I^1},\dots,\omega_{I^{s}}\rangle_{d}\ =\ h^0(\mathcal{M}(I^1,\dots,I^s,d),\Theta)$$
where $\mathcal{M}=\mathcal{M}(I^1,\dots,I^s,d)$ is  a moduli space of parabolic bundles and $\Theta$ an ample line bundle on it. One way to interpret our results is to say that
$$\mathcal{M}\neq \emptyset \Leftrightarrow \langle\omega_{I^1},\dots,\omega_{I^{s}}\rangle_{d}\ =\ h^0(\mathcal{M},\Theta)\neq 0.$$

But $\Theta$ is ample so  $\mathcal{M}\neq \emptyset$ is  equivalent to the assertion $h^0(\mathcal{M},\Theta^N)\neq 0$ for large enough $N$. Hence we obtain
the following generalization of the saturation theorem: For any positive integer $N$, 
\begin{equation}\label{equivalence}
h^0(\mathcal{M},\Theta)\neq 0\ \Leftrightarrow\  h^0(\mathcal{M},\Theta^N)\neq 0.
\end{equation}
In ~\cite{bQ2}, we show that each $h^0(\mathcal{M},\Theta^N)$ is a Gromov-Witten number in a natural way, thereby putting the equivalence ~\ref{equivalence} in a geometric framework.

\section{Properties of  ES-bundles}
\begin{lemma}\label{semi} Let $T$ be a scheme and  $\mw$ a $(D,n)$-vector bundle on $\pone\times\  T$. Suppose we are given a point $t_0\in T$ such that the vector bundle $\mw_{t_0}$ on $\Bbb{P}^1$ is ES.
Then, there exists an open subset  $U\subseteq T$ containing $t_0$ and an  isomorphism $\phi:\mw\ \to\  p_{\pone}^{*}\mz_{D,n}$ over $U$.
\end{lemma}
\begin{proof}
From the hypothesis, there is an isomorphism $s_0:\mw_{t_0}\leto{\sim} \mz_{D,n}$. Now consider $\mt=\shom(\mw,p_{\pone}^*(\mz_{D,n}))$. Clearly $\Hone(\pone,\mt_{t_0})=0$ and therefore (using ~\cite{hart}, Theorem 12.11) $s_0\in H^0(\pone,\mt'_{t_0})$ extends to a neighborhood $U$ of $t_0$:  a  map $\phi:\mw\ \to\  p_{\pone}^{*}\mz_{D,n}$ over $U$ which restricts to  the isomorphism $s_0$ on the fiber over $t_0$. We just need to shrink $U$ further to make $\phi$ into an isomorphism.
\end{proof}
See ~\cite{potier}, Lemma 8.5.3 for the proof of the following lemma:
\begin{lemma}(Serre)\label{serre} Let $\mw$ be  a rank $n$ vector bundle on a smooth projective curve $C$. Assume that $\mw$ is generated by global sections. Then, there exists an exact sequence of the form $$0\ \to\ \mathcal{O}^{n-1}\ \to\  \mw\ \to\  \det(\mw)\ \to\  0.$$
\end{lemma}

\begin{lemma}\label{eric} Let $\lambda$, $D$ and $n$ be integers. There exists an irreducible smooth variety $T$, a vector bundle $\mt$ on $\pone\times\  T$, such that
\begin{enumerate}
\item If $\mw$ is a $(D,n)$-vector bundle on $\pone$ and $\mw=\sum_{i=1}^n \mathcal{O}(a_i)$ with $a_i\geq \lambda$ for $i=1,\dots,n$, then there is a point $t\in T$ and an isomorphism $\mt_t\leto{\sim} \mw$.
\item There is a  Zariski dense open subset $U_{ES}\subseteq T$ formed by points $t$
for which $\mt_t$ is isomorphic to $\mz_{D,n}$.
\end{enumerate}
\end{lemma}
\begin{proof}
Without loss of generality assume (by twisting by an appropriate $\mathcal{O}(a)$) that $\lambda=0$ and hence $D\leq 0$. Let $T=\Ext^1(\mathcal{O}(-D),\oplus_{i=1}^n\mathcal{O})$ and $\mt$ the universal extension. The properties hold because of Lemma~\ref{serre}.
\end{proof}
\begin{lemma}\label{dim1}
 Let $\mv \subseteq \mz_{D,n}$ be a  coherent subsheaf. The following vanishings
 hold
\begin{enumerate}
\item $\Ext^1(\mv,\mz_{D,n})=\Hone(\Bbb{P}^1,\shom(\mv,\mz_{D,n}))=0.$
\item $\ext^1(\mv,\mz_{D,n}/\mv)=\Hone(\Bbb{P}^1,\shom(\mv,\mz_{D,n}/\mv))=0$.
\end{enumerate}
\end{lemma}
\begin{proof} Note the general fact that $\Ext^1(\mf,\mg)=\Hone(\Bbb{P}^1, \shom(\mf,\mg))$ if $\mf$ is a locally free coherent sheaf on $\Bbb{P}^1$ and $\mg$ any coherent sheaf on $\Bbb{P}^1$(see ~\cite{hart}, III, Propositions $6.3$ and $6.7$).

Note that $\mv$ is locally free. Let $\mh=\mz_{D,n}$. By tensoring with a suitable $\mathcal{O}(a)$  we can assume  $\mh=\  \mathcal{O}_{\Bbb{P}^1}^k \oplus\mathcal{O}_{\Bbb{P}^1}(-1)^{n-k}$ for some $k>0$. Let $\mv=\oplus_{i=1}^r\ \mathcal{O}_{\Bbb{P}^1}(a_i)$. It is easy to see that for $i=1,\dots,r$,  $a_i \leq 0$. Hence $\shom(\mv,\mh)$ has a decomposition $\oplus_{u=1}^{rn}\  \mathcal{O}(b_u)$ with each $b_u\geq -1$. But $\Hone(\pone,\mathcal{O}(b))=0$ for $b \geq -1$. This proves $(1)$. Consider the exact sequence of sheaves:
$$0\ \to\ \mv\ \to\ \mh\ \to\ \mh/\mv\ \to\  0$$
which gives an exact sequence
$$0\ \to\ \shom(\mv,\mv)\ \to\ \shom(\mv,\mh)\ \to\ \shom(\mv,\mh/\mv)\ \to\  0$$
 $\Hone(\pone,\shom (\mathcal{V},\mathcal{H}/\mathcal{V}))$ is therefore surjected on by $\Hone(\pone,\shom(\mathcal{V},\mathcal{H}))$ (there are no $H^2$'s on a curve!). This proves $(2)$.
\end{proof}
\subsection{Proof of Proposition ~\ref{bigone}}\label{proofbigone}
 Consider a point  $\mv\subseteq \mh=\mz_{D,n}$ of the moduli space $\mfm(d,r,D,n)$. From Lemma ~\ref{dim1}, $\ext^1(\mv,\mh/\mv)=0$. Therefore (see ~\cite{potier}, Theorem 8.2.1), $\mfm(d,r,D,n)$ is smooth at $\mv$ and the irreducible component of $\mfm(d,r,D,n)$ passing through $\mv$ is of dimension $\chi(\shom(\mv,\mh/\mv))$ which equals $\cee(d,r,D-d,n-r)$.

Consider the subset $U_{G}$ of $\mfm(d,r,D,n)$
formed by  $\mv\subseteq \mh$ such that the bundle $\mv$ is ES. We are going to
show that $U_G$ is open, irreducible and Zariski dense in $\mfm(d,r,D,n)$. This will prove the irreducibility of $\mfm(d,r,D,n)$.

  The openness of $U_G$ follows from  Lemma ~\ref{semi}. For the irreducibility of $U_G$ we argue as follows: Let $Z=Hom(\mz_{d,r},\mz_{D,n})$. On $\pone\times\  Z$, we have a universal morphism
$$\phi:(\mz_{d,r})_Z\ \to\  (\mz_{D,n})_Z.$$
Let $U$ be the largest open subset of $\pone\times\  Z$ such that cokernel of $\phi$ is locally free of rank $n-r$ over $U$. Let $U_Z=Z-p_Z((\pone\times\  Z)-U)$ which is an open set. If a point $z\in Z$ is in $U_Z$, then  $\phi_z$ is an injective map of sheaves on $\pone$ with a locally free cokernel of rank $n-r$. It is easy to see that there is a natural morphism $U_Z \ \to\  \mfm(d,r,D,n)$ and that $U_G$ is the image of this map. Therefore $U_G$ is irreducible.

We now show the Zariski density of $U_G$. Let $\mv\subseteq \mz_{D,n}$ with $\mv=\oplus_{l=1}^r \mathcal{O}(a_l)$ and let $\lambda=\inf_l a_l-1$. We appeal to Lemma ~\ref{eric} for this data and obtain a variety $T$ and a vector bundle $\mt$ on $\pone\times\ T$, and a subset $U_{ES}\subset T$.

 Consider the vector bundle $\ma=\shom(\mt,p_{\Bbb{P}^1}^{*}(\mz_{D,n}))$ on $\pone\times\ T$. Find a $t\in T$ such that $\mt_t$ is isomorphic to $\mv$ and find a section $s_0$ of $H^0(\pone,\ma_t)$ corresponding to the inclusion $\mv\subseteq\mz_{D,n}$. According to Lemma ~\ref{dim1}, $\Hone(\pone,\ma_t)=0$. By Grauert's theorem (see~\cite{hart}, Corollary 12.9) there is an open neighborhood $U$ of $t$ and a section $s$ of $\ma$ on $p_T^{*}(U)$ which restricts to $s_0$ at $t_0$. Said in a different way, we find a map $f:\mt\ \to\  p_{\Bbb{P}^1}^{*}(\mz_{D,n})$ on $\pone\times\ U$ which restricts to the given embedding $\mv\subseteq \mz_{D,n}$ on $\pone\times\ t$. It follows that we can shrink $U$ and assume that
$f$ is injective with locally free cokernel. We therefore have a morphism $\Theta:\ U\ \to\  \mfm(d,r,D,n)$ by the universal property of Quot schemes. Now pick a $t'\in U\cap U_{ES}$  and restrict $f$ to $\pone \times\ t'$. This gives an embedding $\mz_{d,r}\subseteq\mz_{D,n}$. $\Theta(U)$ therefore meets  $U_{G}$ and contains $\mv\subseteq\mz_{D,n}$. Since $U$ is irreducible we have shown the density of $U_G$ in $\mfm(d,r,D,n)$.

By duality, and irreducibility of $\mfm(d,r,D,n)$, the subset  of $\mfm(d,r,D,n)$ formed by $\mv\subseteq\mz_{D,n}$ such that both $\mv$ and $\mz_{D,n}/\mv$ are ES, is open and dense. The proof is therefore complete.
\begin{corollary}\label{forgotten}
The following are equivalent
\begin{enumerate}
\item $\mfm(d,r,D,n)$ is nonempty.
\item $\Hone(\pone,\shom(\mz_{d,r},\mz_{D-d,n-r}))=0.$
\end{enumerate}
\end{corollary}
\begin{proof}
If $\mfm(d,r,D,n)$ is nonempty, pick $\mv\in \mfm(d,r,\mz_{D,n})$ such that both $\mv$ and $\mz_{D,n}/\mv$ are ES. We can now apply Lemma ~\ref{dim1} and deduce $(1)\Rightarrow (2)$.

For $(2)\Rightarrow(1)$, let $\mv=\mz_{d,r}$ and $\mq=\mz_{D-d,n-r}$. Let $\lambda$ be such that $(\mv\oplus\mq)\tensor\mathcal{O}(-\lambda)$ is globally generated. Apply  Lemma ~\ref{eric} to $\lambda, D$ and $n$ and obtain a variety $T$ and a vector bundle $\mt$ on $\pone\times\ T$. By construction, there is a  point $t_{sp}$ on $T$ such that $\mv\oplus\mq=\mt_{t_{sp}}$. We  form the relative Quot scheme $q:\Quot(\mt/T)\ \to\  T$ of quotients of degree $-(D-d)$ and rank $n-r$ of $\mt$. The point $\mv\ \subseteq\  \mv\oplus\mq$ defines a point of $\Quot(\mv\oplus\mq)= \Quot(\mt/T)_{t_{sp}}$. By assumption, $\Hone(\pone,\shom(\mv,\mq))=0$, hence the  Quot scheme $\Quot(\mv\oplus\mq)$ is smooth at $\mv$ of dimension $\chi(\pone,\shom(\mz_{d,r},\mz_{D-d,n-r}))$. Hence we can apply ~\cite{kollar}, Theorem 5.17 to show that $q$ is flat at $\mv$ and hence dominant.  We can now conclude the proof because   $T$ has a nonempty open subset of points $t$ so that  the bundle $\mt_t$ on!
  $\pone$ parameterized by $t$ is ES.
\end{proof}

\section{Proof of (A)$\Rightarrow$ (B) in Theorem ~\ref{MainTe}}\label{easy}
\begin{definition}\label{new101}
Let $\mw$ be a $(D,n)$-bundle on $\pone$, $\me\in \Fl(\mw)$ and $\mi=(d,r,D,n,I)$ a Schubert state. For $p\in\mpp$ let $\pi_p: \mfm(d,r,\mw)\to Gr(d,\mw_p)$ be the natural map. Define the {\em scheme theoretic intersection}
$$\Omega^o(\mi,\mw,\mf)=\bigcap_{p\in\mpp}\pi_p^{-1}[\Omega^o_{I^{p}}(F^p_{\bull})]\subseteq\mfm(d,r,\mw)$$
\end{definition}
It is clear that if $\mw$ is ES and $\mf$ generic, then by  Kleiman's theorem, the dimension of each irreducible component of $\Omega^o(\mi,\mw,\mf)$ is  $\dimii$.

Let $\mv=\mz_{d,r}$ be an ES-bundle and $\mk=(\td,\tr,d,r,K)$ a Schubert state. An element $\mf\in\Fl(\mv)$
is said to be {\bf generic for purposes for intersection theory of $(\mv,\mk)$} if
 the intersection $\Omega^o(\mk,\mv,\mf)$ is proper, and nonempty if and only if $\mk$ is not null. The following Lemma will be proved in Section ~\ref{proofmoving}.
\begin{lemma}\label{moving}\label{openn}
Let $\mw=\mg_{D,n}$, $\mi=(d,r,D,n,I)$ a non null Schubert state, and $\mk=(\td,\tr,d,r,K)$ a non null Schubert state. For $\me$ a generic element of $\Fl(\mw)$, there is an open dense subset of $\Omega^o(\mi,\mw,\me)$
such that for $\mv$ in this open dense subset,
\begin{enumerate}
\item[(i)] $\mv$ and $\mw/\mv$ are ES-bundles.
\item[(ii)] The induced $\me(\mv)\in \Fl(\mv)$ is generic for the purposes of intersection theory of $(\mv,\mk)$.
\end{enumerate}
\end{lemma}
Suppose $\mi=(d,r,D,n,I)$ is a non null Schubert state with $0<r<n$, and $\mk$ a non null Schubert state of the form $(\td,\tr,d,r,K)$ with $0<\tr<r$. We show that Inequality $(\dagger_{\mk}^{\mi})$ from Theorem ~\ref{MainTe} holds:

 Let $\me\in\Fl(\mw)$ be a generic point. Using Lemma ~\ref{moving}, we find
a point in the intersection $\Omega^o(\mi,\mw,\me)$
 which satisfies conditions (i) and (ii) of Lemma ~\ref{moving}.

For $p\in \mpp$ define $L^{p}=\{i^p_a\mid a\in K^{p}\}$ and consider the Schubert state $\ml=(\td,\tr,D,n,L)$. By an easy computation (see ~\cite{ful2}, Lemma 2 (i)),
\begin{equation}\label{biggy}
\Omega^o(\mk,\mv,\me(\mv))\subseteq \Omega^o(\ml,\mw,\me)
\end{equation}
under the inclusion $\mfm(\td,\tr,\mv)\subseteq\mfm(\td,\tr,\mw)$.

Since  $\mk$ is non null (and Lemma ~\ref{moving}), the intersection $\Omega^o(\mk,\mv,\me(\mv))$
is non empty and each irreducible component is of  dimension
$$\dim \mk\ =\ \dim\mfm(\td,\tr,d,r)-\sum_{p\in\mpp}\codim(\omega_{K^{p}}).$$

Since $\mf$ is generic, each irreducible component of $\Omega^o(\ml,\mw,\me)$
is of dimension
$$\dim \ml\ =\ \dim\mfm(\td,\tr,D,n)-\sum_{p\in\mpp}\codim(\omega_{L^{p}}).$$

Therefore inclusion ~\ref{biggy} gives an inequality
$$\dim\mfm(\td,\tr,d,r)-\sum_{p\in\mpp}\codim(\omega_{K^{p}})\leq \dim\mfm(\td,\tr,D,n)-\sum_{p\in\mpp}\codim(\omega_{L^{p}}).$$
This gives
\begin{equation}\label{above112}
[\dim\mfm(\td,\tr,d,r)-\dim \mfm (\td,\tr,D,n)]+\{\sum_{p\in\mpp}\sum_{b=1}^{\tr} (n-\tr+
b-i^p_{k^p_b}-(r-\tr+b-k^p_b))\}\leq 0.
\end{equation}
The term in the square brackets is $-\cee(\td,\tr,D-d,n-r)$ and the term in the curly brackets is $\sum_{p\in\mpp}\sum_{a\in K^p}(n-r+a-i^p_a)$. Therefore Inequality ~\ref{above112} rearranges to Inequality $(\dagger_{\mk}^{\mi})$.
\section{Intersection theory and Universal families}
\subsection{Universal families}
Let $\mi=(d,r,D,n,I)$ be a Schubert state and  $\mw$  a $(D,n)$-vector bundle on $\pone$.
For $T$ a scheme over $\kappa$, let $F(\mi,\mw)[T]$ consist of the set of following data:
\begin{enumerate}
\item A $(d,r)$-subbundle $\mv\subseteq \mw_T$ on $\pone\times\  T$.
\item For each $p\in \mpp$ a complete filtration $E^p_{\bull}$ of the bundle $(\mw_T)_p$ by subbundles (equivalently an element $\me\in\Fl(\mw_T)$).
\end{enumerate}
Subject to the following condition: For each $p\in\mpp$, the rank of
$\mv_p\to (\mw_T)_p/ F^p_{\ell}$ is $r-a$ for $i^p_a\leq \ell< i^p_{a+1}$, $a=0,\dots,r$ ($i^p_0=0$ and $i^p_{r+1}=n$).
It is easy to see that $F(\mi,\mw)$ is a contravariant functor:  schemes/$\kappa$
 to (sets).
\begin{proposition}\label{november11}
The functor $F(\mi,\mw)$ is representable by a scheme $\muu(\mi,\mw)$ which is  smooth over $\mfm(d,r,\mw)$ of fiber dimension
$$\sum_{p\in\mpp}[n(n-2)/2-\codim(\omega_{I^{p}})].$$
\end{proposition}
Note that the dimension of the space of complete flags on a vector space of dimension $n$ is $n(n-1)/2$.
\begin{proof}
Consider an element of $F[T]=F(\mi,\mw)[T]$ as above.
For $p\in\mpp$, the bundle $\mv_p$ on $T$ gets a complete induced filtration (by subbundles). The idea is to rewrite the definition of $F$ by including this data: $F[T]$ is the set of data consisting of
\begin{enumerate}
\item A subbundle $\mv$ of $\mw_T$ on $\pone\times\  T$ of rank $r$ and degree $-d$ when restricted to any fiber $\pone\times\  t$.
\item For each $p\in \mpp$ a complete filtration $E^p_{\bull}$ of the bundle $(\mw_T)_p$ by subbundles.
\item For each $p\in \mpp$ a complete filtration $F^p_{\bull}$ of the bundle $\mv_p$ by subbundles.
\end{enumerate}
Subject to the conditions:
\begin{enumerate}
\item[$(a)$] For $p\in\mpp$, $F^p_{a}\subseteq E^p_{i^p_a}$ for $a=1,\dots,r$
\item [$(b)$] $\mv_p/F^p_a\to (\mw_T)_p/ E^p_{\ell}$ is injective with locally free cokernel for $i^p_a\leq \ell< i^p_{a+1}$ for  $a=1,\dots,r$.
\end{enumerate}
If we ignore the very last condition $(b)$, the desired scheme  is the topmost element $\ma$ in a tower of  Grassmann bundles (see Lemma ~\ref{basic1}) over $Fl_{\mpp}(\tmv)$ where $\tmv$ is the universal subbundle on $\mfm(d,r,\mw)$.
The idea is ``choose the flags on $\mw$ after you have decided what the induced ones on $\mv$ should be''.
The  condition $(b)$  makes the representing scheme an open subscheme of $\ma$.
To calculate the fiber dimension over $M(d,r,\mw)$, it is (using Lemma ~\ref{basic1})
$$\{\sum_{p\in\mpp} r(r-1)/2\}+ \sum_{p\in\mpp}[n(n-1)/2 -\sum_{a=1}^r(i^p_{a+1}-i^p_a)a]$$
(the term in the first curly brackets is the fiber dimension of $\Fl(\mv)$ over
$M(d,r,\mw)$)
$$=\sum_{p\in\mpp} n(n-1)/2+ \{\sum_{p\in\mpp}[1+2+\dots+ r-1 -\sum_{a=1}^r(i^p_{a+1}-i^p_a)a]\}$$
The term in the second curly bracket is  ($i^p_{r+1}=n$)
$$-\sum_{p\in\mpp}[nr-i^p_r-i^p_{r-1}-\dots-i^p_1 -0 -1 -2 -\dots -(r-1)]$$
$$-\sum_{p\in\mpp}[n-r+r-i^p_r +n-r+(r-1)-i^p_{r-1}+\dots n-r+1-i^p_1]
=-\sum_{p\in\mpp}\codim(\omega_{I^{p}})$$
\end{proof}
 For $\me\in \Fl(\mw)$, it is easy to see that $\Omega^o(\mi,\mw,\me)$ is the fiber over $\me$ of the morphism
$\muu(\mi,\mw)\to \Fl(\mw)$.
\begin{corollary}\label{calculstar}
Each irreducible component of $\Omega^o(\mi,\mw,\me)$ is of dimension at least $\dimii$.
\end{corollary}
\begin{proof}
Each irreducible component   of  $\mfm(d,r,\mw)$ passing through a point $\mv$ is of dimension at least
$\chi(\shom(\mv,\mw/\mv))$ (see ~\cite{kollar}, Theorem I.5.17). Each irreducible component of $\muu(\mi,\mw)$ passing through $\mv$ is consequently of dimension at least
$$\chi(\shom(\mv,\mw/\mv))+\sum_{p\in\mpp}[n(n-1)/2-\codim(\omega_{I^{p}})].$$
The dimension of $\Fl(\mw)$ is $\sum_{p\in\mpp}n(n-1)/2$, and
$$\chi(\shom(\mv,\mw/\mv))=\cee(d,r,D,n).$$
Therefore the fiber $\Omega^o(\mi,\mw,\me)$ of the morphism $\muu(\mi,\mw)\to \Fl(\mw)$ over $\me$ is of dimension atleast
$$\chi(\shom(\mv,\mw/\mv))+\sum_{p\in\mpp}[n(n-1)/2-\codim(\omega_{I^{p}})]-\dim(\Fl(\mw))$$
$$\ =\ \chi(\shom(\mv,\mw/\mv))-\sum_{p\in\mpp}\codim(\omega_{I^{p}})\ =\ \dimii.$$
\end{proof}
\begin{corollary}\label{intor}
Suppose $\mw$ is ES, $\me\in \Fl(\mw)$ and $\mv\in \Omega^o(\mi,\mw,\me)$ be such that each irreducible component of $\Omega^o(\mi,\mw,\me)$ which contains $\mv$ is of the expected dimension ($=\dimii$). Then $\mi$ is not null.
\end{corollary}
\begin{proof}
From the smoothness of the schemes $\Fl(\mw)$ and $\muu(\mi,\mw)$, the hypothesis and general statements about flatness (see for example ~\cite{mat}, Theorem 23.1), one concludes that $\muu(\mi,\mw)\to\Fl(\mw)$ is flat at $(\mv,\me)$ and hence dominant. Therefore $\mi$ is not null.
\end{proof}
\begin{corollary}\label{extend}Let $\tmv$ be the universal subbundle on $\mfm(d,r,\mw)$ and $\tmq$ the universal quotient. The morphism
$$\gamma:\muu(\mi,\mw)\to \Fl(\tmv)\times_{\mfm(d,r,\mw)}\Fl(\tmq)$$
is smooth and surjective.
\end{corollary}
\begin{proof}
Let $\mv\subseteq\mw$ be a point in $\mfm(d,r,\mw)$ and $\mq=\mw/\mv$. The surjectivity of the map $\gamma$ is equivalent to the following statement: Given $\mf\in\Fl(\mv)$ and $\mg\in\Fl(\mq)$ then there is a $\me\in\Fl(\mw)$ such that
 $\mv\in\Omega^o(\mi,\mw,\me)$, $\mf=\me(\mv)$ and $\mg=\me(\mq)$.

This is really a pointwise statement (which is left to the reader): If $W$ is a $n$-dimensional vector space, $V$ an $r$-dimensional subspace of $W$, $F_{\bull}\in\Flag(V)$, $G_{\bull}\in\Flag(W/V)$ and $I$  a subset of $[n]$ of cardinality $r$; then there exists $E_{\bull}\in \Flag(W)$ so that $V\in\Omega^o_{I}(E_{\bull})$ and the induced flags on $V$ and $W/V$ are $F_{\bull}$ and $G_{\bull}$ respectively.

The smoothness of $\gamma$ (which is not used in this paper) is proved by the same technique as the smoothness of $\muu(\mi,\mw)\to\Fl(\tmv)$ which was proved
in Proposition ~\ref{november11}.
\end{proof}
\subsection{Proof of Lemma ~\ref{moving}}\label{proofmoving}
We now prove Lemma ~\ref{moving} from Section ~\ref{easy}. Use notation from the statement of Lemma ~\ref{moving}. We first fix a non-empty open subset $U(\mk)$ of $\Fl(\mz_{d,r})$ such that this open subset is invariant under automorphisms of
$\mz_{d,r}$ and such that any $\mf\in U(\mk)$
\begin{itemize}
\item If $\mk$ is null then $\Omega^o(\mk,\mv,\mf)=\emptyset$.
\item If $\mk$ is not null then $\Omega^o(\mk,\mv,\mf)$ is a nonempty proper intersection.
\end{itemize}
Let $\tilde{\mv}$ be the universal subbundle on $\mfm(d,r,\mw)$ and $\Fl(\tmv)$ be the flag bundle as defined in  Section ~\ref{old}. From the previous section we have a smooth morphisms
$$\muu(\mi,\mw)\to \Fl(\tmv)\to \mfm(d,r,\mw).$$

Since $\mi$ is not null, $\mfm(d,r,\mw)$ is non-empty. Therefore the subset $U$
of points in $\mfm(d,r,\mw)$ where the subbundle and the quotient bundle are  $ES$ is nonempty, open  and dense(Proposition ~\ref{bigone}). Let $V\subset \Fl(\tmv)$ be the inverse image of $U$ (clearly non-empty).

Consider the subset $V_{0}$ of $V$ formed by points  $(\mv,\mf)$ so that
$''\mf\in U(\mk)''$. This definition makes sense because if  $(\mv,\mf)\in V$  then
$\mv$ is ES and  $U(\mk)$ is invariant under automorphisms of $\mz_{d,r}$. Locally on $V$, the universal subbundle is ES, and hence the openness of $U(\mk)\subset \Fl(\mz_{d,r})$ implies that $V_0\subset V$ is open. For $\mv\in U\subseteq\mfm(d,r,\mw)$, we can find $\mf\in\Fl(\mv)$ so that $(\mv,\mf)\in V_0$ ($U(\mk)$ is non-empty), and therefore $V_0\neq\emptyset$.

Consider the inverse image (which is a dense open subset) $W_{0}$ of $V_{0}$ under the morphism $\muu(\mi,\mw)\to \Fl(\tmv)$. Using Lemma ~\ref{basic11} we see that for generic $\me\in \Fl(\mw)$, $\Omega^o(\mi,\mw,\me)\cap W_0$ is dense in $\Omega^o(\mi,\mw,\me)$, and this proves Lemma ~\ref{moving}.
\begin{lemma}\label{bprince}\label{basic11} Let $f:X\to Y$ be a morphism of irreducible schemes and  $U_X$ a nonempty open subset of $X$. Then, there exists a nonempty open subset $U_Y$ of $Y$ such that for $y\in U_Y$, $U_X \cap f^{-1}(y)$ is dense in $f^{-1}(y)$.
\end{lemma}
\begin{proof}
Assume without loss of generality that $Y$ is reduced. Let $Z=X-U_X$ with the reduced induced structure. Each irreducible component of $Z$ is of dimension strictly less than $\dim(X)$. Use generic flatness to find a nonempty open subset $U_Y$ of Y such that either $f^{-1}(U_Y)$ is empty or,
\begin{enumerate}
\item $p^{-1}(U_Y)\to U_Y$ is flat and surjective.
\item $p^{-1}(U_Y)\cap Z\to U_Y$ is flat.
\end{enumerate}
It is easy to see that this $U_Y$ satisfies the requirements.
\end{proof}
We note the following corollary to Lemma ~\ref{moving}:
\begin{corollary}\label{trini}
Suppose $\mi=(d,r,D,n,I)$ and $\mk=(\td,\tr,d,r,K)$ are non null Schubert states. Let $\ml=(\td,\tr,D,n,L)$ where $L^{p}=\{i^p_a\mid a\in K^{p}\}$ for $p\in\mpp$. Then
$\ml$ is not null.
\end{corollary}
\begin{proof}
Let $\mw=\mg_{D,n}$. Choose a generic $\mf\in\Fl(\mw)$ and an element $\mv\in\Omega^o(\mi,\mw,\mf)$ so that $\mv$ is ES and $(\mv,\mf(\mv))$ is generic
 for intersection theory of $(\mv,\mk)$.
 $\Omega^o(\mk,\mv,\mf(\mv))$ is therefore non empty. Pick a $\ms$ in this intersection. By an easy calculation, we have $\ms\in \Omega^o(\ml,\mw,\mf)$. Therefore $\ml$ is not null.
\end{proof}
\section{Tangent Spaces}
In this section we denote the Zariski tangent space of a scheme $X$ at a point $x\in X$ by $T(X)_x$. See ~\cite{sottile}, Section 2.7 for the following description of the tangent space of a Schubert variety:
\begin{lemma}\label{sottile}Let $I=\{i_1<\dots<i_r\}$  be a subset of $[n]$ of cardinality $r$ and $W$ an vector space of rank $n$. Let $E_{\sssize{\bullet}}$ be a complete flag on $W$ and $V\in\Omega^{o}_{I}(E_{\sssize{\bullet}})$. Let $E_{\bull}(V)$ and $E_{\bull}(W/V)$ denote the induced flags on $V$ and $W/V$ respectively. Then,
$$T(\Omega^o_I(E_{\bull}))_V \subseteq T(\Gr(r,W))_V= \Hom(V,W/V)$$ is given by
$$\{\phi \in \Hom(V,W/V)\mid \phi(E_a(V)) \subseteq E_{i_a -a}(W/V),\ a=1,\dots, r\}.$$
\end{lemma}
\begin{lemma}\label{tangent}
Let $\mw$ be a $(D,n)$-vector bundle on $\pone$. Let
$\mi$ be a Schubert state of the form $(d,r,D,n,I)$ and $\me\in \Fl(\mw )$. If $\mv\in \Omega^o(\mi,\mw,\me)$ and $\mq=\mw/\mv$ we have
$$T(\Omega^o(\mi,\mw,\me))_{\mv}\ =\
\{\phi\in \Hom(\mv,\mq)\mid \phi_p(E^p_a(\mv))\subseteq E^p_{i^p_a-a}(\mq),\
a=1,\dots,r,\   p\in\mpp \}.
$$
\end{lemma}
\begin{proof}
Let $\me=\prod_p E^p_{\bull}\in \Fl(\mw )$.
The tangent space to $\mfm(r,d,\mw)$ at the point corresponding to $\mv$ is
$\Hom(\mv,\mq)$ (see for example ~\cite{potier}, Theorem 8.2.1). The tangent space to $\Omega^o_{I^{p}}(E^p_{\bull})$ in $Gr(r,\mw_p)$ is described by Lemma ~\ref{sottile}. The lemma now follows from the
 scheme-theoretic description  of $\Omega^o(\mi,\mw,\me)$.
\end{proof}
The vector space given in Lemma ~\ref{tangent} motivates    the following definition.
\begin{definition}\label{SSS}Let $\mi=(d,r,D,n,I)$ be a Schubert state, $\mv$ a $(d,r)$-bundle on $\pone$, $\mq$ a $(D-d,n-r)$-bundle on $\pone$, $\mf\in \Fl(\mv)$ and $\mg\in\Fl(\mq)$. Define
$$\Hom_{\mi}(\mv,\mq,\mf,\mg)=\{\phi\in \Hom(\mv,\mq)\mid \phi_p(F^p_a)\subseteq G^p_{i^p_a-a},\
a=1,\dots,r,\   p\in\mpp \}.
$$
\end{definition}
\begin{lemma}\label{wed}
In the above situation,
\begin{enumerate}
\item For each $(\mf,\mg)\in\Fl(\mv)\times\Fl(\mq)$,
$$
\rk\hoo\geq \dimii\ +\ h^1(\pone,\shom(\mv,\mq)).
$$
\item The set of points $(\mf,\mg)\in \Fl(\mv)\times\Fl(\mq)$ for which
equality holds in (1) is open (possibly empty).
\end{enumerate}
\end{lemma}
\begin{proof}
For each $p\in \mpp$, consider
$$\pi_p:\Hom(\mv,\mq)\to \Hom(\mv_p,\mq_p)$$
and define
$$B_p=\{\psi\in \Hom(\mv_p,\mq_p)\mid \psi(F^p_a)\subseteq G^p_{i^p_a-a},\ a=1,\dots,r\}\ \subseteq \Hom(\mv_p,\mq_p).$$
It is easy to see that codimension of the subspace $B_p\subset \Hom(\mv_p,\mq_p)$ is $\codim(\omega_{I^{p}})$ and $\hoo=\cap_{p\in\mpp}\pi_p^{-1}(B_p)$. The rank of $\hoo$ is therefore at least as large as
$$h^0(\pone,\shom(\mv,\mq))\ -\ \sum_{p\in\mpp}\codim(\omega_{I^{p}})$$
$$=\cee(d,r,D-d,n-r)\ +\ h^1(\pone,\shom(\mv,\mq))\ -\ \sum_{p\in\mpp}\codim(\omega_{I^{p}})$$
$$=\ \dim \mi\ +\ h^1(\pone,\shom(\mv,\mq)).$$
This proves (1). (2) follows from semicontinuity.
\end{proof}
\begin{proposition}\label{tspace}
Let $\mi=(d,r,D,n,I)$ be a Schubert state. Consider the following properties:
\begin{enumerate}
\item[$(\alpha)$] $\mi$ is non null.
\item[$(\beta)$] For generic $(\mf,\mg)\in \Fl(\mz_{d,r})\times\Fl(\mz_{D-d,n-r})$, the rank of
the vector space $\Hom_{\mi}(\mz_{d,r},\mz_{D-d,n-r},\mf,\mg)$ is $\dimii$.
\end{enumerate}
 The following implications hold: $(\beta)\Rightarrow(\alpha)$ in any characteristic and in characteristic $0$, $(\alpha)\Rightarrow(\beta)$.
\end{proposition}
\begin{proof}
{\bf{$(\alpha)\Rightarrow (\beta)$ in characteristic $0$}}: Let $\mw=\mz_{D,n}$.
 Using Lemma ~\ref{below}, we find a $\me\in \Fl(\mw )$ such that
\begin{enumerate}
\item[(1)] $\Omega^o(\mi,\mw,\me)$ is a transverse nonempty intersection.
\item[(2)]  $\Omega^o(\mi,\mw,\me)$ has a dense set subset of points $\mv$  such that both $\mv$ and $\mq=\mw/\mv$ are ES.
\end{enumerate}
Pick a $\mv$ as in $(2)$. We have isomorphisms of bundles $\mv\leto{\sim} \mz_{d,r}$ and $\mw/\mv\leto{\sim} \mz_{D-d,n-r}$.
By $(\alpha)$, Lemmas ~\ref{tspace} and  ~\ref{tangent},  $T(\Omega^o(\mi,\mw,\me))_{\mv}$ is of rank $\dimii$. Hence, using Lemma ~\ref{tangent} and Lemma ~\ref{wed} (2), we find that $(\beta)$ holds.

{\bf{$(\beta)\Rightarrow (\alpha)$ in any characteristic}}: Let $(\mf,\mg)$ be as in $(b)$. Using assumption $(\beta)$ and Lemma ~\ref{wed},
$$H^1(\pone,\Hom(\mz_{d,r},\mz_{D-d,n-r}))\ =\ 0$$
and therefore by Corollary ~\ref{forgotten}, $\mfm(d,r,D,n)\neq\emptyset$. Using Proposition ~\ref{bigone}, find a ES-subbundle $\mv\subseteq \mz_{D,n}$ such that the quotient $\mq=\mz_{D,n}/\mv$ is also ES and choose isomorphisms $\mv\leto{\sim} \mz_{d,r}$ and $\mw/\mv\leto{\sim} \mz_{D-d,n-r}$ and use these to identify $(\mv,\mq)$ with $(\mz_{d,r},\mz_{D-d,n-r})$.

Using Lemma ~\ref{extend}, we can find a $\me \in \Fl(\mw )$ such that
\begin{enumerate}
\item $\mv\in \Omega^o(\mi,\mw,\me).$
\item The flags induced by $\me$ on $\mv$ and $\mq$ are $\mf$ and $\mg$ respectively.
\end{enumerate}
Therefore, by Lemma ~\ref{tangent}  $\Omega^o(\mi,\mw,\me)$ is a transverse and hence  proper intersection at $\mv$. By Lemma ~\ref{intor}, $\mi$ is  not null.
\end{proof}

\begin{lemma}\label{below} If $\kappa$ is algebraically closed of characteristic $0$, $\mw=\mz_{D,n}$ an ES-bundle, $\mi=(d,r,D,n,I)$ a non-null Schubert state and  $\me\in \Fl(\mw)$ is a generic point, the following hold
\begin{enumerate}
\item[(1)] $\Omega^o(\mi,\mw,\me)$ is a transverse nonempty intersection.
\item[(2)]  $\Omega^o(\mi,\mw,\me)$ has a dense open subset of points $\mv$  such that both $\mv$ and $\mq=\mw/\mv$ are ES.
\end{enumerate}
\end{lemma}
\begin{proof}
Property (1) follows from Kleiman's transversality theorem applied to the
morphism ($\Gl(\mw_p)$ acts transitively on $\Gr(r,\mw_p)$ for each $p\in\mpp$)
$$\pi:\mfm(d,r,D,n)\to\prod_{p\in\mpp}\Gr(r,\mw_p).$$

Let $U$ be the open subset of $\mfm(d,r,D,n)$ consisting of $\mv\subset\mw$ such that both $\mv$ and $\mw/\mv$ are ES. This is nonempty by Proposition ~\ref{bigone}. By Kleiman's theorem the dimension of intersection of $\Omega^o(\mi,\mw,\me)$ with the complement of $U$ is of dimension less than $\dim \mi$. Therefore $\Omega^o(\mi,\mw,\me)$ has a dense intersection with $U$.
\end{proof}

\section{Bounds and Genericity}
\subsection{}\label{Bounds}
We begin with an  easy bound:
\begin{lemma}\label{parsul}
Let $\mw$ be a vector bundle on $\pone$. There exists an integer $M$ so that $\deg(\mv)<M$ for any coherent subsheaf $\mv$ of $\mw$.
\end{lemma}
\begin{lemma}\label{bound}
Let $D$, $n$ be integers with $n> 0$. Let $A$ be the set of pairs of (possibly null) Schubert states $(\mi,\mk)$  of the form $\mi=(d,r,D,n,I)$ and $\mk=(\td,\tr,d,r,J)$, such that
\begin{enumerate}
\item[(i)] $\mfm(d,r,D,n)\neq \emptyset$, $\mfm(\td,\tr,D,n)\neq \emptyset$.
\item[(ii)] Inequality $(\dagger_{\mk}^{\mi})$ fails.
\end{enumerate}
Then, the set $A$ is finite.
\end{lemma}
\begin{proof}
The quantity appearing in $(\dagger_{\mk}^{\mi})$ is
$$(-\td )(n-r)+ \tr (D-d)-\tr (n-r)+\sum_p\sum_{a\in K^{p}}(n-r+ a-i^p_a).$$

Therefore if inequality $(\dagger_{\mk}^{\mi})$ fails,
$(-\td )(n-r)+ \tr (D-d)> \tr (n-r)$ and hence $(n-r)(-\td) +\tr(-d)>\tr(n-r)-\tr D$. From Lemma ~\ref{parsul} and hypothesis (i), $-\td$ and $-d$ are  bounded above. So for the finiteness statement, we just need to worry about the case $n=r$ in which case the nonemptiness of $\mfm(d,r,D,n)$ implies that $D=d$ and hence, $(\dagger_{\mk}^{\mi})$ holds.
\end{proof}

Let $\mw=\mz_{D,n}$. Let $B(\mw)\subseteq \Fl(\mw)$ be the largest open subset satisfying the following property: Suppose $\me\in B(\mw)\subseteq\Fl(\mw)$, $\mv$  a subbundle of $\mw$ and $\ms$ a subbundle of $\mv$. Let $\mi$ and $\mk$ be Schubert states defined by  $\mv\in \Omega^o(\mi,\mw,\me)$ and $\ms\in \Omega^o(\mk,\mv,\me(\mv))$. Then, inequality $(\dagger_{\mk}^{\mi})$ holds.

\begin{lemma}\label{unc} $B(\mw)\neq \emptyset.$
\end{lemma}
\begin{proof} Let $U$ be a nonempty open subset of $\Fl(\mw)$ such that for $\me\in U$, $\Omega^o(\mj,\mw,\me)$ is a proper intersection for any  $\mj$ in a fixed finite set which will be made clear below.

Let $\mv$, $\ms$, $\mk$ and $\mi$ be as above. Assume $(\dagger_{\mk}^{\mi})$ fails. For $p\in \mpp$ define $L^{p}=\{i^p_a\mid a\in K^{p}\}$ and consider the Schubert state $\ml=(\td,\tr,D,n,L)$. Then,
\begin{equation}\label{inculsion}
\ms\in\Omega^o(\mk,\mv,\me(\mv))\subseteq \Omega^o(\ml,\mw,\me)
\end{equation}
under the inclusion $\mfm(\td,\tr,\mv)\subseteq\mfm(\td,\tr,\mw)$.
Each irreducible component of $\Omega^o(\mk,\mv,\me(\mv))$ is of dimension atleast (see Corollary ~\ref{calculstar})
$$\dim\mk\ =\ \dim\mfm(\td,\tr,d,r)-\sum_{p\in\mpp}\codim(\omega_{K^{p}}).$$

We will now specify the finite set $\mj$: It is the finite set of  $\ml$ as above obtained from pairs $(\mi,\mk)$
in the set $A$ of Lemma ~\ref{bound}. Hence, each irreducible component of $\Omega^o(\ml,\mw,\me)$ is of
dimension
$$\dim \ml\ =\ \dim\mfm(\td,\tr,D,n)-\sum_{p\in\mpp}\codim(\omega_{L^{p}}).$$
Using Inclusion (~\ref{inculsion}), we therefore have $\dim\mk\leq \dim\ml$ which gives Inequality $(\dagger_{\mk}^{\mi})$ (see Inequality ~\ref{above112}) and hence a contradiction.
\end{proof}

\subsection{List of Genericity properties}\label{list}
Recall that if $\phi:\mv\to \mq$ is a morphism of vector bundles on $\pone$, the kernel of $\phi$ is a
subbundle of $\mv$. This is because the image of $\phi$, being a subsheaf of $\mq$  is locally free. The image of
$\phi$ is however not (necessarily) a subbundle of $\mq$.

Let $\mi=(d,r,D,n,I)$ be a Schubert state (possibly null), $\mv=\mz_{d,r}$ and $\mq=\mz_{D-d,n-r}$. Define $\rk(\mi)$ to be the rank of the vector space
$\hoo$ for generic $(\mf,\mg)\in \Fl(\mv)\times\Fl(\mq)$.

Again, let $(\mf,\mg)$  be a generic point in $\Fl(\mv)\times\Fl(\mq)$ and  $\phi$ a generic element in $\hoo$.  Set $\ms=\ker(\phi)$ and let $\mk=(\td,\tr,d,r,K)$ be the Schubert state determined from the requirement $\ms\in\Omega^o(\mk,\mv,\mf)$. In Section ~\ref{GG} we will prove that the rank of $\hoo$ (which is equal to $\rk(\mi)$ as defined above) is given by the following formula
$$
\dim \mi\ +\ \dim\mk\ +\ \{-\cee(\td,\tr,D-d,n-r)+\sum_{p\in\mpp}\sum_{a\in K^p}(n-r+ a-i^p_a)\}
$$
(the term in the curly brackets below is the quantity appearing in $(\dagger_{\mk}^{\mi})$). Furthermore, $\ms$ and $\mv/\ms$ are ES-bundles and the induced point $(\mf(\ms),\mf(\mv/\ms))\in \Fl(\ms)\times\Fl(\mv/\ms)$ is ``generic''. More precisely, in Section ~\ref{GG}, we construct
\begin{enumerate}
\item[(A)]for arbitrary ES-bundles $\mv$ and $\mq$ on $\pone$, an open subset $A(\mv,\mq)\subseteq B(\mv)\times B(\mq)\subseteq  \Fl(\mv)\times\Fl(\mq)$ (see Section ~\ref{Bounds} for the definition of $B(\mv)$ and $B(\mq)$).
\item[(B)] Given a Schubert state of the form $\mi=(d,r,D,n,I)$, a {\bf non-null} Schubert state of the form $\mk(\mi)=(\td(\mi),\tr(\mi),d,r,K(\mi))$
\end{enumerate}
These satisfy the following properties: Let $\mv=\mz_{d,r}$, $\mq=\mz_{D-d,n-r}$ and $(\mf,\mg)\in A(\mv,\mq)$, then the rank of $\hoo$ is given by the formula (the term in the curly bracket is the term appearing in $(\dagger_{\mk(\mi)}^{\mi})$
\begin{equation}\label{dimformula}
\dim \mi+\dim(\mk(\mi))+\{-\cee(\td(\mi),\tr(\mi),D-d,n-r)+\sum_{p\in\mpp}\sum_{a\in K^p(\mi)}(n-r+ a-i^p_a)\}
\end{equation}
(so that rank of $\hoo$ equals $\rk(\mi)$) and for generic $\phi\in \hoo$, if $\ms=\ker(\phi)$, then
\begin{enumerate}
\item $\ms\in\Omega^o(\mk(\mi),\mv,\mf)$,
\item $\ms$  and $\mv/\ms$ are ES bundles on $\pone$.
\item $(\mf(\ms),\mf(\mv/\ms))\in A(\ms,\mv/\ms)$.
\end{enumerate}

\section{The main technical result}\label{proofes}
\begin{theorem}\label{stonger} Consider a $5$ tuple of the form $(\mv,\mq,\mi,\mf,\mg)$ where $\mv$ is an ES $(d,r)$-bundle and $\mq$ an ES $(D-d,n-r)$-bundle on $\pone$, $(\mf,\mg)$ a generic point of $\Fl(\mv)\times\Fl(\mq)$ and $\mi$ a Schubert state of the form $\mi=(d,r,D,n,I)$.

We claim that there is a filtration by vector subbundles
$${\ms}^{(h)}\subsetneq {\ms}^{(h-1)}\subsetneq\dots \subsetneq {\ms}^{(1)}\subsetneq {\ms}^{(0)}=\mv$$
and injections (of coherent  sheaves) from the graded quotients $\eta_u:{{\ms}^{(u)}}/{{\ms}^{(u+1)}}\hookrightarrow \mq$ for $u=0,\dots, h-1$, such that if we define Schubert states $\mk {(u)}=(\td_u,\tr_u,d,r,K {(u)})$ for $u=1,\dots,h$   by the requirement $\ms^{(u)}\in \Omega^o(\mk {(u)},\mv,\mf)$ then,
\begin{enumerate}
\item[(i)] $\mk {(u)}$, $u=1,\dots,h$ are non-null Schubert states and $\dim \mk {(h)}=0$.
\item[(ii)] For $u=0,\dots,h-1$, $p\in\mpp$ and $a=1,\dots,r$,
$$(\eta_u)_p(\ms^{(u)}_p\cap F^p_a)\subseteq G^p_{i^p_a-a}.$$
\item[(iii)] The rank of $\hoo$ is  (the term in the curly brackets is the term appearing in Inequality $(\dagger_{\mk {(h)}}^{\mi})$)
\begin{equation}\label{aftertheparty}
\dim \mi\ +\ \{-\cee(\td_h,\tr_h,D-d,n-r)\ +\ \sum_{p\in\mpp}\sum_{a\in K^p(h)}(n-r+ a-i^p_a)\}.
\end{equation}
\end{enumerate}
\end{theorem}
\begin{remark}\label{rema11}
If $\mv$ and $\mq$ are ES, the theorem is valid for any $(\mf,\mg)\in A(\mv,\mq)$ (see Section ~\ref{list}).
\end{remark}
\subsection{Proof of (B)$\Rightarrow$ (A) in Theorem ~\ref{MainTe}}\label{second}
We show that Theorem ~\ref{stonger} gives (B)$\Rightarrow$(A) in Theorem ~\ref{MainTe}. Assume condition (B). Let $\mv=\mz_{d,r}$, $\mq=\mz_{D-d,n-r}$, and $(\mf,\mg)$ a generic point of $\Fl(\mv)\times\Fl(\mq)$.

Apply Theorem ~\ref{stonger} to the $5$ tuple $(\mv,\mq,\mi,\mf,\mg)$ and use the same notation. By conclusion (i) of Theorem ~\ref{stonger}, $\mk{(h)}$ is a non null Schubert state.   The dimension of the vector space $\hoo$ is (by conclusion (iii) of Theorem ~\ref{stonger})
$$\dimii+ \{-\cee(\td_h,\tr_h,D-d,n-r)+\sum_{p\in\mpp}\sum_{a\in K^p(h)}(n-r+ a-i^p_a)\}.$$

The hypothesis imply that $(\dagger_{\mk{(h)}}^{\mi})$ holds, therefore
the dimension of $\hoo$ is less than or equal to $\dimii$. We conclude the proof using Proposition ~\ref{tspace}, $(\beta)\Rightarrow(\alpha)$ and Lemma ~\ref{wed}.
\section{Proof of Theorem ~\ref{stonger}}\label{third}
The proof is by induction on $r$. Assume that we have proved this result for all values of $r<r_0$ and prove it for $r=r_0$. For the proof start with $r_0=1$.

If the rank of $\hoo$ is $0$, the filtration is just the singleton $\mv$ and $h=0$, so no maps $\eta$ need to be given. Clearly, the condition in (iii) is met.

So assume the rank of $\hoo$ is nonzero. Pick a generic element $\phi\in \hoo$. Let $\ms$ be the kernel of $\phi$, $\td=-\deg(\ms)$, $\tr=\rk(\ms)$ and
\begin{equation}\label{crackofdawn}
\mk=(\td,\tr,d,r,K)
\end{equation}
the non null Schubert state ($\mk$ is the same as $\mk(\mi)$ from  Section ~\ref{list}) defined by $\ms\in\Omega^o(\mk,\mv,\mf)$. Since $\mk$ is not null, $\dimkk\geq 0$.

  If $\dimkk=0$ take $\ms\subsetneq \mv$ to be the filtration and  $\eta_0=\phi:\mv/\ms \hookrightarrow \mq$. This satisfies (ii) because $\phi\in\hoo$. Since $\dimkk=0$, (i) holds. For (iii), according to Section ~\ref{list}, the rank of $\hoo$ is (since $\dimkk=0$),
$$\dimii+ \{-\cee(\td,\tr,D-d,n-r)+\sum_{p\in\mpp}\sum_{a\in K^p}(n-r+ a-i^p_a)\}.$$

Therefore assume that $\dimkk>0$ and therefore  $0<\tr<r$. Now by our discussion of genericity (Section ~\ref{list}) we may apply the induction hypothesis on the $5$-tuple $(\ms,\mv/\ms,\mk,\mf(\ms),\mf(\mv/\ms))$.

 We therefore find a filtration
${\ms}^{(h)}\subsetneq {\ms}^{(h-1)}\subsetneq\dots\subsetneq {\ms}^{(1)}=\ms$
 and morphisms $\gamma_u:{\ms}^{(u)}/{\ms}^{(u+1)}\hookrightarrow {\mv}/{\ms}$ for $u=1,\dots, h-1$ which satisfy the conclusions of Theorem ~\ref{MainTe} for the 5-tuple $(\ms,\mv/\ms,\mk,\mf(\ms),\mf(\mv/\ms))$.

We  claim that the filtration $${\ms}^{(h)}\subsetneq {\ms}^{(h-1)}\subsetneq\dots\subsetneq {\ms}^{(1)}={\ms}\subsetneq {\ms}
^{(0)}=\mv$$ and maps $\eta_u= \phi\gamma_u$ for $u=1,\dots,s$, $\eta_0=\phi$ satisfy the conditions (i), (ii) and (iii) in the theorem.

 For $u=1,\dots,h$, let ${\Ftwo}{(u)}$  the Schubert state defined from $\ms^{(u)}\in\Omega^o(\Ftwo{(u)},\mv,\mf)$ (${\Ftwo}{(1)}$ is same as $\mk$ defined above). Also for $u=1,\dots,h$, let  $\ml{(u)}$ be the Schubert state defined by $\ms^{(u)}\in \Omega^o(\ml{(u)},\ms,\mf(\ms))$.

For $u=1,\dots,h$, let $\td_u$, $\tr_u$, $L(u)$ and $K(u)$ be determined from
 $$\ml(u)=(\td_u,\tr_u,\td,\tr,L(u))$$
 $$\mk(u)=(\td_u,\tr_u,d,r,K(u))$$
(so $\tr_u$ is the rank of ${\ms}^{(u)}$ and  $\td_u=-\deg(\ms^{(u)})$.)

By an easy calculation and Equation ~\ref{crackofdawn},
\begin{equation}\label{night}
K^p(u)=\{k^p_b\mid b\in L^p(u)\}
\end{equation}
{\bf Verification of (i)}: We know that $\mk$ is not null (Section ~\ref{list}) and by the induction hypothesis,  each $\ml(u)$ is not null. Therefore by Corollary ~\ref{trini}, and Equation ~\ref{night}, $\mk(u)$ is not null for $u=0,\dots,h$.

Inductive conclusion (iii) for the five tuple $(\ms,\mv/\ms,\mk,\mf(\ms),\mf(\mv/\ms)$ tells us that the rank of $\hoos$ is
$$\dimkk+ \{-\cee(\td_h,\tr_h,d-\td,r-\tr)+\sum_{p\in\mpp}\sum_{b\in L^p{(h)}}(r-\tr+ b-k^p_b)\}$$
The term in the curly brackets is the quantity appearing in $(\dagger_{\ml(h)}^{\mk})$ which is  $\ \leq 0$. However, by Lemma ~\ref{wed}, the rank of $\hoos$ is
at least as great as $\dimkk$. We therefore conclude that equality holds in Inequality $(\dagger_{\ml(h)}^{\mk})$:
\begin{equation}\label{E1}
-\cee(\td_h,\tr_h,d-\td,r-\tr)+\sum_{p\in\mpp}\sum_{b\in L^p{(h)}}(r-\tr+ b-k^p_b)=0.
\end{equation}
Or that,
\begin{equation}\label{E11}
\cee(\td_h,\tr_h,d-\td,r-\tr)-\sum_{p\in\mpp}\sum_{t=1}^{d_h}(r-\tr+ l^p(h)_t-k^p_{l^p(h)_t})=0.
\end{equation}
Also inductive conclusion (i) for $\mk$ says that
\begin{equation}\label{c1}
\dim(\ml{(h)})=0
\end{equation}
 This is the same as the equation
\begin{equation}\label{E2}
\cee(\td_h,\tr_h,\td-\tilde{d}_h,\tr-\tr_h)-\sum_{p\in\mpp}\sum_{t=1}^{d_h}(\tr-\tr_h+t-l^p(h)_t)
=0.
\end{equation}
Adding Equation ~\ref{E11} to Equation ~\ref{E2} and using
Equation ~\ref{night}, we get the desired conclusion:
\begin{equation}\label{E3}
\dim\ \Ftwo(h)\ =\ \cee(\td_h,\tr_h,d-\td_h,r-\tr_h)-\sum_{p\in\mpp}\sum_{t=1}^{d_h}(r-\tr_h+t
-k^p_{l^p(h)_t})=0.
\end{equation}
{\bf Verification of (ii):} We need to verify that for ${u}=0,\dots,h-1$, $p\in\mpp$ and $a=1,\dots,r$,
\begin{equation}\label{dag}
(\eta_{u})_p({\ms}^{({u})}_p\cap F^p_a)\subseteq G^p_{i^p_a-a}.
\end{equation}
Now suppose ${u}$, $p$ and $a$ are as above. If ${u}=0$, Inclusion ~\ref{dag}
is clear because $\phi\in \hoo$. So assume ${u}>0$ and find $t$ such that
$F^p_a\cap  \ms_p=F^p_t(\ms)$. Clearly, $k^p_t\leq a$.

From $\eta_{{u}}=\phi \gamma_{u}$, we see that
$$(\eta_{{u}})_p({\ms}^{({u})}_p\cap F^p_a)=\phi_p(\gamma_{u})_p({\ms}^{({u})}_p\cap F^p_t(\ms))\subseteq \phi_p (F^p_{k^p_t-t}({\mv}/{\ms}))$$
where in the last inclusion, we have used the property (ii) satisfied by the maps $\gamma_{u}$: $(\gamma_{u})_p({\ms}^{({u})}_p\cap F^p_t(\ms))\subseteq F^p_{k^p_t-t}({\mv}/{\ms}).$ But,
$$\phi_p (F^p_{k^p_t-t}({\mv}/{\ms}))=\phi_p((F^p_{k^p_t}+\ms_p)/{\ms_p})=\phi_p(F^p_{k^p_t})\subseteq \phi_p(F^p_a)\subseteq G^p_{i^p_a-a}.$$
{\bf Verification  of (iii):}
We claim,
\begin{claim}\label{pope}
For $u=1,\dots,h$, let $b(u)$ be the quantity
\begin{equation}\label{defb}
\dim\ \mk(u)-\dim\ \ml(u)\ +\ \{-\cee(\td_u,\tr_u,D-d,n-r)+\sum_{p\in\mpp}\sum_{a\in K^p(u)}(n-r+ a-i^p_a)\}.
\end{equation}
Then $b(u)\leq b(u+1)$ for $u=1,\dots,h-1$.
\end{claim}
Note that the term in the curly brackets in Quantity ~\ref{defb} is the
 same as the one in Inequality ($\dagger_{\mk(u)}^{\mi}$).

The claim implies that $b(h)\geq b(1)$. But $\dim \ml(1) =0$ and therefore by Section ~\ref{list}, $b(1)+\dimii$ is the rank of the vector space $\hoo$. Hence, the rank of the vector space $\hoo$ is not greater than $\dimii+b(h)$. Also, from (i) and induction, $\dim\ \mk(h)=\dim\ \ml(h)=0$. Therefore the rank of $\hoo$ is less than or equal to  Expression ~\ref{aftertheparty}.

But according to Lemma  ~\ref{later}, the rank of $\hoo$ is at least as much as  Expression ~\ref{aftertheparty}. In conjunction with the above, this says that the rank of $\hoo$ is equal to Expression ~\ref{aftertheparty}. The proof of Theorem ~\ref{stonger} would therefore be complete once the claim is proved.

{\bf{Proof of the claim:}}
We would like to (compare with proof of Claim 5.2 in ~\cite{jag}) apply  Lemma ~\ref{unc}, to $\mg\in B(\mq)$ and  $\im(\eta_u)\subseteq\im(\phi)\subseteq \mq$.  However, we cannot do so, because the image of $\eta_u$ (and $\im(\phi)$) may not be a subbundle of $\mq$. Such difficulties are routine in the theory of parabolic bundles where one proves that the saturation of the image subsheaf has ``better properties'' than the subsheaf.

For $u=0,\dots,h$,
let ${\mee}^{(u)}$ be the saturation of the image of $\eta_u$ in $\mq$,  and let $\mee=\mee^{(0)}$ be the saturation of $\text{image} (\phi)$ in $\mq$ (see Section ~\ref{satu} for the notion of saturation of a coherent subsheaf of a vector bundle on a smooth curve).
$$\frac{\ms^{(u)}}{\ms^{(u+1)}}\letoo{\eta_u}\mee^{(u)}\subseteq\mee^{(0)}\subseteq\mq.$$

The claim will be proved by applying  Lemma ~\ref{unc} to $\mg\in B(\mq)$ and  $\mee^{(u)}\subseteq\mee\subseteq \mq$ (see Remark ~\ref{rema11}).

It is important to bring in the concept of parabolic bundles in our calculation (see Section ~\ref{satu} for the notation and basic results). Before resuming the proof in Section ~\ref{illinois} of Claim ~\ref{pope}, we write the  inequalities $(\dagger_{\mk}^{\mi})$ in terms of parabolic degrees in the next section.
\subsection{Parabolic Stability and the Horn problem}\label{ParHor}
Let $\mw$ be an ES-bundle of degree $-D$ and rank $n$, $\me\in\Fl(\mw)$ and $\mv\subseteq\mw$ a subbundle. Let $\mi=(d,r,D,n,I)$ be the  Schubert state determined from the condition $\mv\in \Omega^o(\mi,\mw,\me)$.

Let $\me(\mv)=\prod_{p\in\mpp}E^p_{\bull}(\mv_p)$ be the induced collection of flags on $\mv$. Define weights as follows: Let $w^p_a=\frac{n-r+a-i^p_a}{n-r}$ for $p\in\mpp$ and $a=1,\dots,r$. This gives a structure of a parabolic bundle $\umv=(\mv,\me(\mv),w)$ on $\mv$.

Consider a subbundle $\ms\subseteq\mv$ and  $\mk=(\td,\tr,d,r,K)$ be the Schubert state determined from the requirement $\ms\in\Omega^o(\mk,\mv,\me(\mv))$.
\begin{equation}\label{ustad}
\pa(\ms,\umv)=-(n-r)\td+\sum_{p\in\mpp}\sum_{a\in K^p}(n-r+a-i^p_a)
\end{equation}
$$= [-(n-r)\td+\cee(\td,\tr,D-d,n-r)]+\{-\cee(\td,\tr,D-d,n-r)+\sum_{p\in\mpp}\sum_{a\in K^p}(n-r+a-i^p_a)\}$$
The term in the curly brackets is the same as in $(\dagger_{\mk}^{\mi})$. The term in the square brackets is $\tr(n-r)-(D-d)\tr=\chi(0,\tr,D-d,n-r)$. We conclude that if $\me\in B(\mw)$ (see Lemma ~\ref{unc}), we have
\begin{equation}\label{ustad4}
(n-r)\pa(\ms,\umv)\leq \chi(0,\tr,D-d,n-r)
\end{equation}
We divide by $\tr(n-r)$ and write Equation ~\ref{ustad4} as
 \begin{equation}\label{ustad3}
\mupar(\ms,\umv)\leq \frac{1}{\tr(n-r)}\chi(0,\tr,D-d,n-r)=1 -\frac{D-d}{n-r}.
\end{equation}
\subsection{Return to the proof of Claim ~\ref{pope}}\label{illinois}
Let $\underline{\mee}$ the induced parabolic bundle  corresponding (see Section ~\ref{ParHor})  to the subbundle $\mee\subseteq\mq$ and the given collection of flags $\mg\in \Fl(\mq)$.

 Suppose that the rank of $\ker(\phi_p)\supseteq\ms_p$ is $\tr +\epsilon(p)$ (recall that the rank of $\ms$ is $\tr$) and  $H^p=\{h^p_1<\dots<h^p_{\tr+\epsilon(p)}\}$  the unique subset of $\{1,\dots,r\}$ such that 
\begin{equation}\label{RelativePosition}
\ker(\phi_p)\in \Omega^o_{H^{p}}(F^p_{\bull})\subseteq\Gr(\tilde{r}+\epsilon(p),\mv_p).
\end{equation}
Introduce the notation
$$\theta^p(a)=i^p_a-a$$ for $p\in\mpp$, $a=1,\dots,r$,
 $$\bed=D-d,\ \emm=n-r,$$
$$\gammar=r-\tr,\  \degr=d-\td.$$
$$\dm=-\deg(\mee),$$
$$c\ =\ \frac{1}{n-r-(r-\tr)}\ =\ \frac{1}{\emm-\gammar}.$$

We note the following inequality for the degree of $\mee$:
\begin{equation}\label{v111}
-\dm\geq (-d)-(-\td)+\sum_{p\in\mpp} \epsilon(p)=-\degr+\sum_{p\in\mpp}\epsilon(p).
\end{equation}
\begin{lemma}\label{l101}
For $a\in[\tr+\epsilon(p)]$ and $p\in\mpp$, the weight attached to $\mee_p\cap G^{p}_{\theta^p(h^p(a))}(\mq_p)$ (which is a member of the induced flag $G^p_{\bull}(\mee)$) in the parabolic bundle $\underline{\mee}$ is at least
$$c[{\emm}-\gammar+ (h^p_a-a)-\theta^p(h^p_a)].$$
\end{lemma}
\begin{proof}
To see this, fix a $p$ and suppose  $\mee_p\in\Omega^o_J(G^p_{\bull})$ where $J$ is a subset $J=\{j_1<\dots<j_{\rk(\mee)}\subset \{1,\dots,{\emm}\}$. Now, if $\rk(\mee_p\cap G^{p}_{\theta^p(h^p_a)}(\mq_p))=x$ and $x\neq 0$ then $j_x\leq \theta^p(h^p_a)$. Also, $x\geq h^p_a-a$  because
\begin{enumerate}
\item $\phi_p(F^p_{h^p_a}(\mv_p))$ is $h^p_a-a$ dimensional.
\item $\phi_p(F^p_{h^p_a}(\mv_p))\subseteq \mee_p\cap G^{p}_{\theta^p(h^p_a)}(\mq_p).$
\end{enumerate}
If $x=0$ which could happen only when $h^p_a=a$,
$$c[{\emm}-\gammar+ (h^p_a-a)-\theta^p(h^p_a)]\leq c({\emm}-\gammar)\leq 1$$
 and the weight $w^p_0\geq 1$.
\end{proof}
\begin{lemma}\label{l102}
The parabolic degree $\pa({\mee}^{(u)},\underline{\mee})$, is  at least
$$(\td_{u+1}-\td_u)+c\sum_{p\in\mpp}\sum_{{\temp}\in L^p(u)-L^p(u+1)}[{\emm}-\gammar+ (k^p_{\temp}-{\temp})-\theta^p(k^p_{\temp})]-c\sum_{p\in\mpp}\epsilon(p)(\tr_u-\tr_{u+1}).$$
\end{lemma}
\begin{proof}
For $p\in\mpp$,  suppose that $$\ms_p\in\Omega^o_{U^{p}}(F^p_{\bull}(\ker(\phi_p))),\ U^{p}=\{u^p_1<\dots<u^p_{\tr}\}\subseteq[\tilde{r}+\epsilon(p)].$$

We have (using Equation ~\ref{RelativePosition}),  $K^p=\{h^p_{u^p_{\temp}}\mid {\temp}=1,\dots,\tr\}$ for $p\in\mpp$. Therefore $\phi_p$ maps $F^p_{k^p_{\temp}}(\mv_p)$ to an element of the flag  $G^p_{\bull}(\mee_p)$ whose weight is at least (by Lemma ~\ref{l101})
$$c[{\emm}-\gammar+ (h^p_{u^p_{\temp}}-u^p_{\temp})-\theta^p(h^p_{u^p_{\temp}})]$$
$$=c[{\emm}-\gammar+ k^p_{\temp}-u^p_{\temp}-\theta^p(k^p_{\temp})].$$
$$=c[{\emm}-\gammar+ k^p_{\temp}-{\temp}-\theta^p(k^p_{\temp})]+c[{\temp}-u^p_{\temp}].$$
But clearly ${\temp}-u^p_{\temp}\geq -\epsilon(p).$ Therefore the above quantity is
$$\geq c[{\emm}-\gammar+ (k^p_{\temp}-{\temp})-\theta^p(k^p_{\temp})]-c[\epsilon(p)].$$

For ${\temp}=1,\dots,\tr$, by (ii), $(\eta_u)_p(F^p(\ms)_{\temp}\cap {\ms_p}^{(u)})\subseteq \mee_p\cap G^{p}_{\theta^p(k^p_{\temp})}$. Hence by Lemma ~\ref{naive}, the statement follows.
\end{proof}
We now return to the proof of the claim. The fundamental inequality that we are going to use is Equation ~\ref{ustad4} (since $\mg\in B(\mq)$ by assumption). Recall that $\mee$ is a $(\dm,{\gammar})$- bundle.
\begin{equation}\label{eq2}
({\emm}-\gammar)\pa({\mee}^{(u)},\underline{\mee})-\chi(0,\tr_u-\tr_{u+1},\bed-\dm,{\emm}-\gammar)\leq 0.
\end{equation}
Using Lemma ~\ref{l102} we obtain the inequality (use $c({\emm}-\gammar)=1$),
$$(\emm-\gammar)(\td_{u+1}-\td_u)+\sum_{p\in\mpp}\sum_{{\temp}\in L^p(u)-L^p(u+1)}[{\emm}-\gammar+ (k^p_{\temp}-{\temp})-\theta^p(k^p_{\temp})]$$
\begin{equation}\label{july18}-\sum_{p\in\mpp}\epsilon(p)(\tr_u-\tr_{u+1}) -\chi(0,\tr_u-\tr_{u+1},\bed-\dm,{\emm}-\gammar)\leq 0.
\end{equation}

The first, third and fourth term in the above inequality combine to give
$$({\emm}-\gammar)(\td_{u+1}-\td_u)-\chi(0,\tr_u-\tr_{u+1},\bed-\dm,{\emm}-\gammar)-\sum_{p\in\mpp}\epsilon(p)(\tr_u-\tr_{u+1})$$
$$=-\chi(\td_{u}-\td_{u+1},\tr_u-\tr_{u+1},\bed-\dm-\sum_{p\in\mpp}\epsilon(p),{\emm}-\gammar)$$
\begin{equation}\label{eq3}
\geq -\chi(\td_u-\td_{u+1},\tr_u-\tr_{u+1},\bed\ {-\ {\degr}},{\emm}-\gammar).
\end{equation}
(using Inequality ~\ref{v111} in the last step)

We now deduce from  Inequality ~\ref{july18} and Inequality ~\ref{eq3}   that
\begin{equation}\label{mainone}
\sum_{p\in\mpp}\sum_{{\temp}\in L^p(u)-L^p(u+1)}[{\emm}-\gammar+ (k^p_{\temp}-{\temp})-\theta^p(k^p_{\temp})]
-\chi(\td_u-\td_{u+1},\tr_u-\tr_{u+1},\bed\ -\ \degr,{\emm}-\gammar)\leq 0.
\end{equation}
Write the left hand side of Inequality ~\ref{mainone} as
$A(u)-A(u+1)$ where
$$A(u)=\sum_{p\in\mpp}\sum_{{\temp}\in L^p(u)}[{\emm}-\gammar+ (k^p_{\temp}-{\temp})-\theta^p(k^p_{\temp})]-\chi(\td_u,\tr_u,\bed-\degr,{\emm}-\gammar)$$$$=B(u)+C(u)$$
with
$$B(u)=\chi(d_u,r_u,\degr,{\gammar})-\sum_{p\in\mpp}\sum_{{\temp}\in L^p(u)}[\gammar-(k^p_{\temp}-{\temp})]=\dim(\mk(u))-\dim(\ml(u)).$$
 (see derivation of Inequality ~\ref{above112}) and
$$C(u)=-\chi(d_u,r_u,\bed,{\emm})+\sum_{p\in\mpp}\sum_{{\temp}\in L^p(u)}[{\emm}-\theta^p(k^p_{\temp})]$$
$$=-\cee(\td_u,\tr_u,D-d,n-r)+\sum_{p\in\mpp}\sum_{a\in K^p(u)}(n-r+ a-i^p_a).$$
 Therefore, $b(u)=A(u)$, we conclude by observing that  Inequality ~\ref{mainone} gives $A(u+1)\geq A(u)$.

\subsection{Proof of Theorem ~\ref{MainTe}}
In Theorem ~\ref{MainTe}, we have already proved $(A)\Rightarrow (B)$ (Section ~\ref{easy}) and  $(B)\Rightarrow(A)$ (Section ~\ref{second}). The implications  $(B)\Rightarrow(C)\Rightarrow(D)$ are obvious and we are left with having to prove $(D)\Rightarrow(B)$. We will assume the basic transversality result (Proposition ~\ref{carol}, which is independent of this section).

Suppose (by way of contradiction) that (D) holds and (B) fails and  $\ml$ is a non null Schubert state for which inequality $(\dagger_{\ml}^{\mi})$ fails.

Let $\mv$ be a ES $(d,r)$-bundle and $\mf$ a generic point of $\Fl(\mv)$. Consider the parabolic structure on $\mv$ where to $F^p_a$ we assign the weight
$w^p_a=\frac{n-r+a-i^p_a}{n-r}$. This gives us a parabolic bundle $\umv$.

Let $\mk=(\td,\tr,d,r,K)$ be  a non null Schubert state. Let $\ms\in\Omega^o(\mk,\mv,\mf)$ (which is nonempty by genericity of $\mf$).

We have the following formula for the parabolic slope (see Equation ~\ref{ustad}). The term in the curly brackets is the same as the one in inequality $(\dagger_{\mk}^{\mi})$
\begin{equation}\label{u2}
\tr(n-r)\mupar(\ms,\umv)= [\tr(n-r)-(D-d)\tr]+\{-\cee(\td,\tr,D-d,n-r)+\sum_{p\in\mpp}\sum_{a\in K^p}(n-r+a-i^p_a)\}.
\end{equation}
 The assumption $\dimii\geq 0$ gives us that (take $\mk$ to
correspond to $\ms=\mv$ above) $\mupar(\mv,\umv)\leq 1-\frac{D-d}{n-r}$. Let $\mt\in \Omega^o(\ml,\mv,\mf)$, then
by Equation ~\ref{u2},
$$\mupar(\mt,\umv)>1-\frac{D-d}{n-r}\geq\mupar(\mv,\umv).$$
Therefore $\umv$ is not a semistable parabolic bundle. Let $\mathcal{X}$ be the
Harder-Narasimhan maximal contradictor of semistability (see ~\cite{MS}). It has the following
properties: If $\mathcal{J}$ is the Schubert state determined from $\mathcal{X}\in \Omega^o(\mj,\mv,\mf)$, then $\langle \mj\rangle =1$
and $\mupar(\mathcal{X},\umv)$ is the maximum slope among all subbundles of $\mv$ (and among the ones that reach this maximum, it has the maximum rank). The  equality $\langle \mj\rangle =1$ is consequence of the uniqueness of $\mathcal{X}$ and the genericity of $\mf$ (see ~\cite{b1}) and the transversality result (Proposition ~\ref{carol}).

Therefore \begin{equation}\label{funny}
\mupar(\mathcal{X},\umv)\geq\mupar(\mt,\umv)>1-\frac{D-d}{n-r}.
\end{equation}
But $\langle \mj \rangle =1$ and therefore by assumption (D), inequality $(\dagger_{\mj}^{\mi})$ holds.
Applying Equation ~\ref{u2}, we therefore find that $\mupar(\mathcal{X},\umv)\leq 1-\frac{D-d}{n-r}$ which is in contradiction to Equation ~\ref{funny}.
\section{Some additional remarks}
For the purpose of future use, we would like to note the following result which follows from the proof of Theorem ~\ref{stonger}. This section can be skipped in the first reading (it is not used in this paper).

Suppose that $\mv$ is an {\bf arbitrary} $(d,r)$-bundle, $\mq$ an ES $(D-d,n-r)$-bundle on $\pone$ and $\mi$ a Schubert state of the form $\mi=(d,r,D,n,I)$. Let $\mf\in \Fl(\mv)$ be {\bf arbitrary} and $\mg\in B(\mq)$ (or just ``generic'').

Assume that there is a sequence
$${\ms}^{(h)}\subsetneq {\ms}^{(h-1)}\subsetneq\dots \subsetneq {\ms}^{(1)}\subsetneq {\ms}^{(0)}=\mv$$
and injections (of coherent  sheaves) from the graded quotients $\eta_u:{{\ms}^{(u)}}/{{\ms}^{(u+1)}}\hookrightarrow \mq$ for $u=0,\dots, h-1$ (such a sequence can be obtained inductively as in Theorem ~\ref{stonger}). Assume further that for $u=0,\dots,h-1$, $p\in\mpp$ and $a=1,\dots,r$,
$$(\eta_u)_p(\ms^{(u)}_p\cap F^p_a)\subseteq G^p_{i^p_a-a}.$$

Define Schubert states $\mk {(u)}=(\td_u,\tr_u,d,r,K {(u)})$ and $\ml(u)$
for $u=1,\dots,h$   by the requirement $\ms^{(u)}\in \Omega^o(\mk {(u)},\mv,\mf)$ and $\ms^{(u)}\in \Omega^o(\ml{(u)},\ms,\mf(\ms))$.

Then the conclusions of Claim ~\ref{pope} holds:
\begin{theorem}\label{stoner}
For $u=1,\dots,h$ let
\begin{equation}\label{defbb}
b(u)=\dim(\mk(u))-\dim(\ml(u))+\{-\cee(\td_u,\tr_u,D-d,n-r)+\sum_{p\in\mpp}\sum_{a\in K^p(u)}(n-r+ a-i^p_a)\}.
\end{equation}
Then, $b(u)\leq b(u+1)$ for $u=1,\dots,h-1$.
\end{theorem}
 In the proof of Claim ~\ref{pope},we only used that $\mg\in B(\mq)$, the genericity of $\mf$ played no role. Theorem ~\ref{stoner} follows immediately.
\section{Dimension counts and genericity}\label{GG}
\subsection{Rank of $\hoo$}\label{trinitie}

\begin{lemma}\label{rain}
Let $\mv$ and $\mq$ be bundles on $\pone$. There exist only finitely many pairs $(\td,\tr)\in\Bbb{Z}^2$ such that
there exist a $(\td,\tr)$ subbundle $\ms\subseteq\mv$ and an injection of sheaves $\pi:\mv/\ms\to \mq$.
\end{lemma}
\begin{proof}
By Lemma ~\ref{parsul}, degrees of $\ms$ and $\mv/\ms$ are both bounded above. Therefore the set of possible degrees of $\ms$ is both bounded above and below.
\end{proof}

Use notation from Definition ~\ref{SSS}, and assume that $\mv$ and $\mq$ are ES-bundles. Suppose $\phi$ is a generic element of $\hoo$ and $\ms=\ker(\phi)$. It is easy to see that $\ms$ is a {\em subbundle} of $\mv$. Suppose $\ms\in \Omega^o(\mk,\mv,\mf)$ for a Schubert state $\mk=(\td,\tr,d,r,K)$.

By Lemma ~\ref{rain} are only finitely many  possible degrees and ranks for  the kernel of $\phi$ (given $\mv$ and $\mq$).
Lemma ~\ref{rain} also implies that $\mk$ is not null. This is because by the genericity of $\mf$, we may assume that each of the intersections $\Omega(\mh,\mv,\mf)$ is empty if $\mh$ is a null schubert state of the form $(d_1,r_1,d,r,K)$ satisfying
\begin{itemize}
\item There is a homomorphism  $\phi:\mv\to \mq$ such that the kernel of $\phi$ is of degree $-d_1$ and rank $r_1$.
\end{itemize}
(There are only finitely many such $\mh$ because of Lemma ~\ref{rain}.)

 The aim of this section is to prove the following proposition:
\begin{proposition}\label{june}
With notation and assumptions as above, $\hoo$ is of rank
\begin{equation}\label{firstformula}
\dim\mi\ +\ \dim\mk\  +\  \{\sum_{p\in\mpp} \sum_{a\in K^{p}}(n-r+a-i^p_a)\  -\ \cee(\td,\tr,D-d,n-r)\}.
\end{equation}
where the term in the last curly brackets is the quantity in Inequality $(\dagger_{\mk}^{\mi})$.
\end{proposition}
Introduce for convenience the notation,
$$\bedda=D-d,\ \emman=n-r.$$
 Let $\Hom(\mv,\mq,\td,\tr)$ be the space of morphisms
$\mv\to \mq$ such that the kernel is a subbundle of rank $\tr$ and degree $-\td$ (the scheme structure will be given in Section ~\ref{trin}). We postpone the proof of the following Lemma to Section ~\ref{trin}.
\begin{lemma}\label{one}
$\Hom(\mv,\mq,\td,\tr)$ is smooth and irreducible (possibly empty) of
rank
$$\dim\mfm(\td,\tr,\mv)\ +\ \cee(d-\td,r-\tr,{\bedda},{\emm})$$
\end{lemma}
 For each $p\in\mpp$, let $Z(p)$ be the smooth  scheme with points $(\ess,\psi)$ where
\begin{itemize}
\item $\ess$ is a $\tr$ dimensional subspace of $\mv_p$.
\item $\psi$ is a morphism $\mv_p/\ess\to \mq_p$ (not necessarily injective).
\end{itemize}
There is a map  $\eta_p:\Hom(\mv,\mq,\td,\tr)\to Z(p)$ given by $\phi\leadsto((\ker(\phi))_p,\phi_p)$.  Let $T(p)\subseteq Z(p)$ be the subscheme of points $(\ess,\psi)$ such that
\begin{enumerate}
\item For $a=1,\dots,r$, $\psi(F^p_a)\subseteq G^p_{i^p_a-a}$.
\item $\ess\in \Omega^o_{K^{p}}(F^p_{\bull})$.
\end{enumerate}
It is easy to see that $T(p)\to\Omega^o_{K^{p}}(F^p_{\bull})$ is a smooth morphism of fiber dimension $\sum_{a\in [r]\smallsetminus K^{p}}(i^p_a-a).$ It now follows that $T(p)$ is smooth and
\begin{equation}\label{az12}
\codim(T(p),Z(p))=\codim(\omega_{K^{p}})+(r-\tr)\emm-\sum_{a\in [r]\smallsetminus K^{p}}(i^p_a-a).
\end{equation}
Let $\Omega=\bigcap_{p\in\mpp}\eta_p^{-1}(T(p)).$ It is easy to see that $\Omega$ is a dense open subset of $\hoo$. Each irreducible component of $\Omega$ is of dimension at least

$$\rk\Hom(\mv,\mq,\td,\tr)-\sum_{p} \codim(T(p),Z(p))=$$
$$\dim\mfm(\td,\tr,d,r)+\cee(d-\td,r-\tr,{\bedda},{\emm})-\sum_{p\in\mpp}[\codim(\omega_{K^{p}})+((r-\tr)\emm-\sum_{a\in [r]\smallsetminus K^{p}}(i^p_a-a))]=$$

$$\{\dim\mfm(\td,\tr,d,r)-\sum_{p\in\mpp}\codim(\omega_{K^{p}})\}+\{\cee(d-\td,r-\tr,{\bedda},{\emm})-\sum_{p\in\mpp}((r-\tr)\emm +\sum_{a\in [r]\smallsetminus K^{p}}(i^p_a-a))\}$$
The term in the first curly bracket is $\dim\mk$  and the second term is
$$\cee(d,r,{\bedda},{\emm})-\sum_{p\in\mpp} \sum_{a\in [r]}({\emm}+a-i^p_a)
-\cee(\td,\tr,{\bedda},{\emm})+\sum_{p\in\mpp} \sum_{a\in K^{p}}({\emm}+a-i^p_a).$$
$$=\dimii+\{-\cee(\td,\tr,{\bedda},{\emm})+\sum_{p\in\mpp} \sum_{a\in K^{p}}({\emm}+a-i^p_a)\}.$$
{\em Therefore the rank of $\hoo$ is at least as much the quantity ~\ref{firstformula}.}

However, the above only gives us an inequality. Kleiman's theorem cannot be applied because the schemes $Z(p)$ may not be homogenous for the group $\Gl(\mv_p)\times \Gl(\mq_p)$.

But one can get an exact expression for the rank of $\hoo$ if track is kept
of the kernels of the morphism $\phi_p:(\mv/\ms)_p\to \mq$. For $p\in\mpp$,
\begin{itemize}
\item Let $\epsilon(p)$ be the rank of the kernel $\mb^p\subset (\mv/\ms)_p$ of the map $\phi_p:(\mv/\ms)_p\to \mq$.
\item Define $J^{p}\subset[r-\tr]$ a subset of cardinality $\epsilon(p)$ from the requirement $\mb^p\in \Omega^o_{J^{p}}(F^p_{\bull}(\mv/\ms))$.
\end{itemize}
Let $\Hom(\mv,\mq,\td,\tr,\epsilon)$ be the space of morphisms $\psi:\mv\to \mq$ such that the kernel is a subbundle of rank $\tr$ and degree $-\td$ and such that for $p\in\mpp$, the kernel of $\psi_p:\mv_p\to\mq_p$ is of rank $\tr+\epsilon(p)$. We postpone the proof of the following lemma to Section ~\ref{trin}.
\begin{lemma}\label{two}
$\Hom(\mv,\mq,\td,\tr,\epsilon)$ is smooth and irreducible of
dimension
$$\dim\mfm(\td,\tr,\mv)+\cee(d-\td,r-\tr,{\bedda},{\emm})-\sum_{p\in\mpp}\epsilon(p)(\emm-(r-\tr-\epsilon(p)))$$
\end{lemma}
Consider $\ma=\Hom(\mv,\mq,\td,\tr,\epsilon)$. For each $p\in\mpp$
, let $Y(p)$ be the scheme with points $(\ess,B,\psi)$ where
\begin{itemize}
\item $\ess$ is a $\tr$ dimensional subspace of $\mv_p$.
\item $B$ is a $\epsilon(p)$ dimensional subspace of $\mv_p/\ess$.
\item $\psi$ is a morphism $\mv_p/\ess\to\mq_p$ with kernel $B$.
\end{itemize}
There is a morphism $\lambda_p:\ma\to Y(p)$ given by $$\phi\leadsto((\ker(\phi))_p, ker(\phi_p)/(\ker(\phi))_p,\phi_p).$$
Clearly, the group $GL(\mv_p)\times GL(\mq_p)$ acts transitively on $Y(p)$. Now fix a $p\in\mpp$ and assume that we are given complete flags $F_{\bull}(\mv_p)$ and $G_{\bull}(\mq_p)$ on $\mv_p$ and $\mq_p$ respectively. Let $R(p)\subseteq Y(p)$ be the subscheme of points $(\ess,B,\phi)$ such that
\begin{enumerate}
\item $\ess\in \Omega^o_{K^{p}}(F^p_{\bull})$.
\item $B\in \Omega^o_{J^{p}}(F^p_{\bull}(\mv_p/\ess)).$
\item For $a=1,\dots,r$, $\psi(F^p_a)\subseteq G^p_{i^p_a-a}$ (where we consider $\psi$ as a morphism $\mv_p\to\mq_p$).
\end{enumerate}
For $p\in\mpp$, let
$$[r]\smallsetminus K^{p}=\{\alpha^p(1)<\dots<\alpha^p(r-\tr)\}$$
 It is easy to see that the codimension of $R(p)$ in $Y(p)$ is (fiber $R(p)$ over the set of choices of $(C,B)$):
$$\codim(\omega_{K^{p}})+\codim(\omega_{J^{p}})+[(r-\tr-\epsilon(p))\emm-\sum_{t\in [r-\tr]\smallsetminus J^{p}}(i^p_{\alpha^p(t)}-\alpha^p(t))]$$
and
$$\sum_{t\in[r-\tr]\smallsetminus J^{p}}(i^p_{\alpha^p(t)}-\alpha^p(t))=\sum_{t\in[r-\tr]}(i^p_{\alpha^p(t)}-\alpha^p(t))-\sum_{t\in J^{p}}(i^p_{\alpha^p(t)}-\alpha^p(t)).$$
Also note that
\begin{equation}\label{star11}
\sum_{t\in[r-\tr]}\ttt=\sum_{\ell\in [r]\smallsetminus K^{p}}(i^p_{\ell}-{\ell}).
\end{equation}
Therefore  $\codim(R(p),Y(p))$ is
\begin{equation}\label{return}
\sum_{\ell\in[r]\smallsetminus K^p}(\emm +\ell- i^p_{\ell})-\sum_{t\in J^{p}}(i^p_{\alpha^p(t)}-\alpha^p(t))-\epsilon(p)\emm+\codim(\omega_{K^{p}})+\codim(\omega_{J^{p}})
\end{equation}
It follows from our assumptions that
$$\bigcap_{p\in\mpp}\lambda_p^{-1}(R(p))\subset\hoo$$
for generic $\mf$ and $\mg$ is a dense subset of dimension (by Kleiman's theorem)
$\dim(\ma)-\sum_p(\codim(R(p),Y(p)))$. Use the expression for $\dim(\ma)$ given  in Lemma ~\ref{two} and Equation ~\ref{return} to write this as a sum of two parts. The first one  is
$$\dim(\mfm(\td,\tr,\mv))-\sum_{p\in\mpp}\codim(\omega_{K^{p}})+\cee(d-\td,r-\tr,{\bedda},{\emm})-\sum_{p\in\mpp}\sum_{\ell\in[r]\smallsetminus K^p}[\emm+\ell-i^p_{\ell}]$$
and the second part is sum over $p\in\mpp$ of the ``local discrepancies''
$$\Disc(p)=-\epsilon(p)({\emm}-(r-\tr-\epsilon(p))-\codim(\omega_{J^{p}})+\epsilon(p)\emm-\sum_{t\in J^{p}}\ttt.$$
The first part is the integer from the Expression ~\ref{firstformula}. Fix a $p$ and use the notation $J^{p}=\{j_1<\dots<j_{\epsilon(p)}\}$. Then,
$$\codim(\omega_{J^{p}})=\sum_{a=1}^{\epsilon(p)}[(r-\tr)-\epsilon(p)+a-j_a].$$
Therefore the local ``discrepancy'' at $p$ is
$$\Disc(p)=\sum_{a=1}^{\epsilon(p)}(-\emm+(r-\tr-\epsilon(p)) -(r-\tr-\epsilon(p)+a-j_a)+{\emm}
-(i^p_{\alpha^p(j_a)}-\alpha^p(j_a)))$$
$$=\sum_{a=1}^{\epsilon(p)}(j_a-a+\alpha^p(j_a)-i^p_{\alpha^p(j_a))}).$$
 Pick $\phi\in \bigcap\lambda_p^{-1}(R(p))\subset\hoo$. For each $t$,
$$\phi_p (F^p_t(\mv/\ms)) \subseteq G^p_{\ttt}(\mq).$$
Take $t=j_a$, the rank of the left hand side is $j_a-a$ and the right hand side  is of rank $i^p_{\alpha^p(j_a)}-\alpha^p(j_a)$. Therefore $$j_a-a\leq i^p_{\alpha^p(j_a)}-\alpha^p(j_a)$$ and hence $\Disc(p)\leq 0$.

We therefore obtain that the rank of $\hoo$ for generic $(\mf,\mg)$ is less than or equal to Expression ~\ref{firstformula}. Hence, for generic $(\mf,\mg)$ the rank of $\hoo$  is given by Expression ~\ref{firstformula} (it is not clear who to thank for the success of this argument!). This concludes the proof of Proposition ~\ref{june}.
\subsection {A rank inequality}
The following inequality was used in Section ~\ref{third}:
\begin{lemma}\label{later}
Suppose that  $(\mv,\mq,\mf,\mg,\mi)$ is  as in Definition ~\ref{SSS}. Suppose $\mk$ is a Schubert state of the form $(\td,\tr,d,r,K)$ and $\Omega^o(\mk,\mv,\mf)\neq \emptyset.$ Then the rank of $\hoo$ is at least
\begin{equation}\label{exppress}
\dimii+\{-\cee(\td,\tr,{\bedd},{\emma})+\sum_{p\in\mpp}\sum_{\ell\in K^p}({\emma}+ \ell-i^p_{\ell})\}.
\end{equation}
\end{lemma}
\begin{proof}
Pick  $\ms\in \Omega^o(\mk,\mv,\mf)$. The rank of $\Hom(\mv/\ms,\mq)$ is at least  $$\chi(\shom(\mv/\ms,\mq))= \cee(d-\td,r-\tr,{\bedd},{\emma}).$$
This inequality is strict if $\Hone(\pone,\shom(\mv/\ms,\mq))\neq 0$.

Consider,
$\hos=\{\phi\in \hoo\mid \phi(\ms)=0\}.$ For each $p\in \mpp$ we have a map $\eta_p :\Hom(\mv/\ms,\mq)\to \Hom(\mv_p/\ms_p,\mq_p)$. Let $W^p=\{\eta\in \Hom(\mv_p/\ms_p,\mq_p)\mid \phi_p(F^p_a)\subseteq G^p_{i^p_a-a},\  a=1,\dots,r\}$. The rank of $W^p$ is $\sum_{a\in [r]\smallsetminus K^{p}}(i^p_a-a)$ and hence its codimension in $\Hom(\mv_p/\ms_p,\mq_p)$ is $(r-\tr)({\emma})-\sum_{a\in [r]\smallsetminus K^{p}}(i^p_a-a)=\sum_{a\in [r]\smallsetminus K^{p}}({\emma}+a-i^p_a)$. It is easy to see that $\hos$ is the intersection of vector spaces $\bigcap\eta_p^{-1}W^p$ and is consequently of rank at least
$$
\cee(d-\td,r-\tr,{\bedd},{\emma})-\sum_{p\in\mpp}\sum_{a\in [r]\smallsetminus K^{p}}({\emma}+a -i^p_a)$$
$$=\cee(d,r,{\bedd},{\emma})-\sum_{p\in\mpp}\sum_{a=1}^r({\emma}+a-i^p_a)-\cee(\td,\tr,{\bedd},{\emma})+\sum_{p\in\mpp}\sum_{\ell\in K^{p}}({\emma}+\ell-i^p_{\ell})$$
which is the same as Expression ~\ref{exppress}.
\end{proof}
\subsection{Proofs of Lemma ~\ref{one} and ~\ref{two}}\label{trin}
 Let $\mv$ be a $(d,r)$-bundle and $\mq$ a $({\bedd},{\emma})$-bundle on $\pone$. Use once again the notation $\bedda=D-d,\emman=n-r$. Let
 $\epsilon:\mpp\to \Bbb{Z}_{\geq 0}$, and $\td$, $\tr\in \Bbb{Z}$.

 We define $F(\mv,\mq,\td,\tr,\epsilon)$ to be the contravariant functor:  schemes/$\kappa\to$ (sets) which assigns to every scheme $T$ over $\Spec(\kappa)$ the set of data of the form $(\ms,\mb,\phi)$ of the form
\begin{enumerate}
\item[(a)]$\ms$ is a $(\td,\tr)$-subbundle of $\mv_T$.
\item[(b)] An injection of sheaves $\phi:\mv_T/\ms\to \mq_T$ such that the quotient of $\phi$ is flat over $T$.
\item[(c)] $\mb=\prod_{p\in\mpp} \mb^p$ where for $p\in\mpp$,  $\mb^p$ is a subbundle of $(\mv_T/\ms)_p$ of rank $\epsilon(p)$.
\end{enumerate}
The above data is required to satisfy: For $p\in\mpp,\ \ker(\phi_p)=\mb^p$ and $\phi_p:(\mv_T/\ms)_p\to (\mq_T)_p$ has a locally free cokernel (of rank $\rk(\mq)-\rk(\mv)+\rk(\ms)+\epsilon(p)$).

Denote the functor obtained by considering pairs $(\ms,\phi)$ satisfying (a) and (b) above by $F(\mv,\mq,\td,\tr)$. We make the following observation:
\begin{lemma}\label{friday}
Suppose  that $\mq$ is an  ES-bundle. If $T$ is a scheme and  $(\ms,\phi)$ satisfies $(a)$ and $(b)$ above then $$\Hone(\pone,\shom(\mv/\ms_t,\mq))=0,\  \forall\  t\in T.$$
\end{lemma}
\begin{proof}
Via $\phi$, $\mv/\ms_t\subseteq\mq$ is an sheaf theoretic injection of vector bundles. But $\mq$ is an $ES$ bundle, so we can apply Lemma ~\ref{dim1}.
\end{proof}
Consider $H=Hom(\mv,\mq)$ which represents the functor:  schemes/$\kappa \to$ (sets) given by $T\leadsto Hom(\mv_T,\mq_T)$ (see Lemma ~\ref{somelemma}). We have a universal morphism $\phi:\mv_H\to\mq_H$ on $H$. let $\mt$ be the cokernel of this morphism.
By the method of flattening stratifications (see ~\cite{curves}, Lecture 8) applied to $\mt$ and each $\mt_p$ for $p\in\mpp$, it is clear that each $F(\mv,\mq,\td,\tr,\epsilon)$ (and $F(\mv,\mq,\td,\tr)$) is representable by a scheme which we denote by $\Hom(\mv,\mq,\td,\tr,\epsilon)$ (and respectively $\Hom(\mv,\mq,\td,\tr)$). These schemes coincide on the level of points with the sets of the same name in Section ~\ref{trinitie}.

 To motivate the functorial arguments in this section, let us first parameterize the points of
$\Hom(\mv,\mq,\td,\tr)$ when both $\mv$ and $\mq$ are ES-bundles. That is, in the definition above take
$T=\Spec{\kappa}$. We first select $\ms$: The set of choices is parameterized by $\mfm(\td,\tr,\mv)$.

We have to now parameterize the possible choices $A(\ms)$ of $\phi:\mv/\ms\to \mq$ so that $(\ms,\phi)$ satisfy
the requirements. We ``write down'' all morphisms $\mv/\ms\to \mq$. This is a vector space
$H(\ms)=\Hom(\mv/\ms,\mq)$. We would not be interested in $\ms$ if $H(\ms)$ did not contain a sheaf-theoretic
injection. If it did, the rank of $H(\ms)$ is $\chi(\shom(\mv/\ms,\mq))$ because the $H^1$ vanishes (Lemma
~\ref{dim1}) and $A(\ms)$ is the open subset of $H(\ms)$ formed by injections.
We will show that the  set of $\ms$ such that $A(\ms)$ contains an injection is an open subset of $\mfm(\td,\tr,\mv)$ (essentially because of the vanishing $H^1$).

To parameterize $\Hom(\mv,\mq,\td,\tr,\epsilon)$ when both $\mv$ and $\mq$ are ES bundles, we proceed similarly:
we select $\ms$ and $\mb\in\prod_{p\in\mpp} \Gr(\epsilon(p),(\mv/\ms)_p)$. We then shift $\mv/\ms$ along $\mb$
(see Section ~\ref{shiftop}), and proceed as before (consider morphisms of the shift $\Psi(\mv/\ms,\mb)\to\mq$).

Lemmas ~\ref{one} and ~\ref{two} follow from
\begin{lemma}\label{represent}\label{thanksgiving}\label{cor}Assume that $\mv$, $\mq$, $d$, $r$ and $\epsilon$ are as above and that in addition that $\mv$ and $\mq$  are ES-bundles. Let $\my$ be the universal quotient on $\pone\times \mfm(\td,\tr,\mv)$.
Then,
\begin{enumerate}
 \item $\Hom(\mv,\mq,\td,\tr)$ and  $\Hom(\mv,\mq,\td,\tr,\epsilon)$ are smooth  and irreducible schemes over $\kappa$.
\item Consider the diagram
\begin{equation}\label{dig}
\xymatrix{\Hom(\mv,\mq,\tr,\td,\epsilon)\ar@{^{(}->}[r]^{i}\ar[d]^{a}&  \Hom(\mv,\mq,\tr,\td)\ar[ddl]^{\pi}\\
 \prod_p Gr(\epsilon(p),\my_p)\ar[d]^{t}\\
\mfm(\td,\tr,\mv)\\}
\end{equation}
The morphisms $\pi$, $a$ and $t$ are each smooth (therefore all objects above are smooth
 and irreducible) and $i$ is an inclusion.
Also,
\begin{enumerate}
\item $\Hom(\mv,\mq,\tr,\td,\epsilon)$ is irreducible (possibly empty) of dimension
$$\dim\mfm(\td,\tr,\mv)+\cee(d-\td,r-\tr,{\bedda},{\emman})-\sum_{p\in\mpp} \epsilon(p)[\emman-(r-\tr-\epsilon(p))].$$
\item  $\Hom(\mv,\mq,\tr,\td)$ is irreducible (possibly empty) of dimension
$$\dim\mfm(\td,\tr,\mv)+\cee(d-\td,r-\tr,{\bedda},{\emman}). $$
\end{enumerate}
\end{enumerate}
\end{lemma}
\begin{proof}
Let $X=\mfm(\td,\tr,\mv)$. Let $\my$ be the universal quotient of $\mv_X$ on $\pone\times\  X$. Consider the vector bundle $\mt=\shom(\my,\mq_X)$ on $\pone\times\  X$.
 Let $U_X$ be the open subset of $X$ formed by points $x$ such that $\Hone(\pone,\mt_x)=0$. Let $Y$ be the total space of $(p_T)_*\mt$ (over $U_X$). It is clear from the definition of $F(\mv,\mq,\td,\tr)$ and Lemma ~\ref{somelemma}, Lemma ~\ref{friday}, that there is a natural morphism
 $\Hom(\mv,\mq,\td,\tr)\to Y.$

 On $\pone\times\ Y$ we have a morphism of sheaves $$\phi:p_{\pone\times X}^*(\my)\to \mq.$$ Let $U_Y$ be the open subset of points $y$ of $Y$ such that $\phi_y$ is injective as a morphism of sheaves (Lemma ~\ref{aa}). It is easy to see that  $U_Y$ represents $F(\mv,\mq,\td,\tr)$.

To represent $F(\mv,\mq,\td,\tr,\epsilon)$ we form the Grassmann bundle
$$\tilde{X}=\prod_{p\in\mpp}Gr(\epsilon(p),\my_p)$$
over $X$. Over $\tilde{X}$, for each $p\in\mpp$ we have bundles $\mb^p\subseteq \my_p$. We can therefore form the shift (a bundle on $\pone\times\  \tilde{X}$, see Section ~\ref{shiftop}): $\tilde{\my}=\Psi(\mb,\my)$. Let $\tilde{\mt}=\shom(\tilde{\my},p_{\pone}^*(\mq))$ on $\pone\times\ \tilde{X}$. Let $U_{\tilde{X}}$ be the open subset of $\tilde{X}$ formed by $\tilde{x}$ such that
$\Hone(\pone,\tilde{\mt}_{\tilde{x}})=0$. Let $\tilde{Y}$ be the total space of $p_{\tilde{X} *}\tilde{\mt}$ over $U_{\tilde{X}}$. On $\pone\times\tilde{Y},$ there is a morphism
$$\tilde{\phi}:p_{\pone\times \tilde{X}}^*(\tilde{\my})\to \mq.$$ Let $U_{\tilde{Y}}$ be the open subset of points $\tilde{y}$ of $\tilde{Y}$ such that $\tilde{\phi}_{\tilde{y}}$ is injective as a morphism of sheaves (see Lemma ~\ref{injective} for the openness of this set).
 We find on $\pone\times U_{\tilde{Y}}$ an injective (over each point of $U_{\tilde{Y}}$) morphism of sheaves
$$p_{\pone\times\tilde{X}}^*(\tilde{\my})\to \mq.$$
For each $p\in \mpp$, consider the (composite) morphism on $U_{\tilde{Y}}$
$$\phi_p:p_{\tilde{X}}^*(\my_p) \to p_{\tilde{X}}^*(\tilde{\my}_p)\leto{\tilde{\phi}_p} \mq_p.$$
This morphism has $p_{\tilde{X}}^*(\mb^p)$ in its kernel. We look at the open subset of $U_{\tilde{Y}}$ formed by points where the cokernel of $\phi_p$ is locally free and the kernel is exactly $p_{\tilde{X}}^*(\mb^p)$ for each $p\in \mpp$. By Lemma ~\ref{somelemma} and Section ~\ref{patrick}, this open subset  represents $F(\mv,\mq,\td,\tr,\epsilon)$.

The smoothness assertions are clear from the above arguments.
For the first dimension formula (the second one is obtained similarly): From the proof, we need to
add the Euler characteristic of $\shom(\tilde{\my},\mq)$ to the dimension of the
Grassmann bundle $\tilde{X}$. The degree of $\tilde{\my}$ is $-(d-\td)+\sum_p\epsilon(p)$. The Euler characteristic of $\shom(\tilde{\my},\mq)$
is  therefore $\cee(d-\td,r-\tr,{\bedda},{\emman}) - \sum_{p\in\mpp}\epsilon(p)\emman$
and the dimension of $\tilde{X}$ is
$$\dim\mfm(\td,\tr,\mv)+\sum_{p\in\mpp}\epsilon(p)(r-\tr-\epsilon(p)).$$
\end{proof}
\begin{remark}
 In Lemma ~\ref{cor} it is easy to see that $t$ is surjective. However $a$ is a smooth dominant morphism which need not be surjective. This is because the map $U_{\tilde{X}}\to \Gr(\td,\tr,\mv)$ may not be surjective.
\end{remark}
\subsection{Open Subsets in Products of Flag Spaces}\label{on1}
We now construct the open subsets in products of flag spaces as required by Section ~\ref{list}.

Let $\mv$, $\mq$, $\mi$ and $\mk$, $d$, $r$, $\epsilon$ and $H^p$ (for $p\in\mpp$) be as in section ~\ref{trinitie}.  Form a universal scheme $\muu_{gen}(\mv,\mq,\mi)$ as a scheme over $\Fl(\mv)\times \Fl(\mq)$ whose fiber over a point $(\mf,\mg)$ is the subset of $\Hom(\mv,\mq,\td,\tr,\epsilon)$ formed by $\phi$ that if we set $\ms=\ker(\phi)$ and for $p\in\mpp$, $\mb^p\subseteq(\mv/\ms)_p $ the kernel of  $\phi_p:(\mv/\ms)_p\to \mq$, then
\begin{itemize}
\item $\ms\in\Omega^o(\mk,\mv,\mf)$.
\item For $p\in\mpp$, $\epsilon(p)=\rk(\mb^p)$ and $\mb^p\in \Omega^o_{H^{p}}(F^p_{\bull}(\mv/\ms))$.
\end{itemize}
We prove the following lemma in the next section.
\begin{lemma}\label{dominant1}
The natural morphism $\mug\to\muu(\mk,\mv)$ is smooth and dominant. $\mug$ is a smooth scheme (recall that we are assuming that $\mv$ and $\mq$ are $ES$ bundles).
\end{lemma}
\begin{remark}
The possible non surjectivity of $\mug\to \muu(\mk,\mv)$ is related the morphism $a$ from Diagram ~\ref{dig} being non surjective.
\end{remark}
Let $O(\mi)\subset\Fl(\mv)\times\Fl(\mq)$ be the largest open subset such that if $(\mf,\mg)\in O(\mi)$, the dimension count of Section ~\ref{trinitie} applies.

We now construct the desired open subset $A(\mv,\mq)$ from section ~\ref{list}. The Schubert state $\mk(\mi)$ required is just $\mk$ from above. It is natural to require functoriality under isomorphisms for $A(\mv,\mq)$. That is, $A(\mv,\mq)$ should be stable under the automorphisms of $\mv$ and $\mq$.

If $r=1$ define $A(\mv,\mq)=\Fl(\mv)^s\times B(\mq)$. Assume that the construction has been achieved if $r<r_0$ and we extend it to the case $r=r_0$. The requirements are over a finite set of choices of $\mi$. It is enough to satisfy the requirements for each $\mi$, obtain open subsets $A(\mv,\mq,\mi)$ and finally set
$$A(\mv,\mq)=\bigcap_{\mi}A(\mv,\mq,\mi).$$

Fix a $\mi$. Consider the diagram
$$\xymatrix{\mug\ar[d]^{p}\ar[r]\ar[r]^{\pi} & \Fl(\mv)\times\Fl(\mq)\\
            \muu(\mk,\mv)\\
            }$$
Here $p$ is dominant (Lemma ~\ref{dominant1}). Let $U$ be the largest open subset of $\muu(\mk,\mv)$ formed by points $(\ms,\mf)$ such that
\begin{enumerate}
\item $\ms$ and $\mv/\ms$ are ES (Section ~\ref{moving}).
\item $(\mf(\ms),\mf(\mv/\ms))\in A(\ms,\mv/\ms)$ (Lemma ~\ref{extend}).
 \end{enumerate}

Shrink $O(\mi)$ so that $\mug\to \Fl(\mv)\times\Fl(\mq)$ is flat over it. Find the largest possible such $O(\mi)$.

We finally
define $$A(\mv,\mq,\mi)=\pi(p^{-1}(U))\cap O(\mi)\cap[B(\mv)\times B(\mq)]$$
which is clearly an open and nonempty subset of $\Fl(\mv)\times\Fl(\mq)$ and satisfies the requirements.
\subsection{Proof of Lemma ~\ref{dominant1}}
Let $\mi$, $\mk$, $\epsilon$, $d$, $r$ be as in Section ~\ref{on1}. We want to show that the natural morphism $p:\mug\to\muu(\mk,\mv)$ is smooth (and hence dominant because $\mug\neq \emptyset$).

Define $\mathcal{C}$ by the cartesian square
$$\xymatrix
{\mathcal{C}\ar[r]\ar[d] & \muu(\mk,\mv)\ar[d]\\
\hoq\ar[r] &  \mfm(\td,\tr,\mv)}$$

It is easy to see that $\mathcal{C}$ parameterizes triples $(\ms,\phi,\mf)$ so that $(\ms,\mf)\in\muu(\mk,\mv)$,
$\ker(\phi)=\ms$ and $\phi\in\hoq$. The morphism $\hoq\to\mfm(\td,\tr,\mv)$ is smooth and dominant (it is the
composition $t\circ a$  in Diagram ~\ref{dig} of Lemma ~\ref{represent}). Therefore $\mathcal{C}\to
\muu(\mk,\mv)$ is smooth and dominant.

Therefore for Lemma ~\ref{dominant1}, we need to show that $\mug\to\mathcal{C}$ is a smooth morphism. In fact it
is (Zariski locally) a fiber bundle with smooth fibers. We will argue ``at the level of points'' (to be rigorous
one work at the level of fibers of functors as in previous section, but we leave this to the reader). What is the
fiber of $\mug\to\mathcal{C}$ over $(\ms,\phi,\mf)$? It is just the set of
$\mg=\prod_{p\in\mpp}G^p_{\bull}\in\Fl(\mq)$ such that $\forall p\in\mpp$, $a\in[r]$, $\phi_p(F^p_{a})\subseteq
G^p_{i^p_a-a}.$ The rank of $\phi_p(F^p_a)$ can be expressed by a formula in terms of the given data (independent
of $\phi$). Therefore using Lemma ~\ref{basic1} we conclude that the set of $\mg$ is parameterized by a smooth
space whose dimension doesn't jump as we vary $(\ms,\phi,\mf)$. This concludes the proof of Lemma
~\ref{dominant1}.

\section{Transversality}\label{gamma}
The following proposition was obtained in some special cases by F. Sottile ~\cite{sottile},~\cite{sottile2}. The classical part of it was proven independently by the author ~\cite{jag} and R. Vakil ~\cite{vakil}.
\begin{proposition}\label{carol}
Let $\mw$ be an ES $(D,n)$-bundle on $\pone$ and $\mi$ a Schubert state of the form $\mi=(d,r,D,n,I)$.  For generic $\mf\in\Fl(\mw)$, the smooth points of $\Omega^o(\mi,\mw,\mf)$ are dense in $\Omega^o(\mi,\mw,\mf)$. This implies that if $\dimii=0$,  $\Omega^o(\mi,\mw,\mf)$ is a smooth and transverse intersection of dimension $0$.
\end{proposition}
\begin{proof}
If $\mi$ is null, there is nothing to prove. Therefore assume that $\mi$ is not null.

We first show that there is a dense open subset $\hat{U}(\mi,\mw)$ of $\muu(\mi,\mw)$ such that for $(\mv,\me)\in \hat{U}(\mi,\mw)$, $\Omega^o(\mi,\mw,\me)$ is a transverse intersection at $\mv$.

 Pick a point $\mv\in\mfm(d,r,\mw)$ (which is nonempty because $\mi$ is non null) such that
$\mv$ and $\mq=\mw/\mv$ are both ES. Pick a generic point $(\mf,\mg)\in\Fl(\mv)\times\Fl(\mq)$. Now find (using Lemma ~\ref{extend}) $\me\in\Fl(\mw)$ such that $\mv\in\Omega^o(\mi,\mw,\me)$ and the collection of flags induced on $\mv$ and $\mq$ are $\mf$ and $\mg$ respectively.

Since, $\mi$ is not null, condition (B) of Theorem ~\ref{MainTe} holds for
 the Schubert state $\mi$. By the proof of $(B)\Rightarrow(A)$ in Section ~\ref{second}, the tangent space to $\Omega^o(\mi,\mw,\me)$ at $\mv$ which is given by $\hoo$ is of rank $\dimii$. We have therefore found a point $(\mv,\me)$ such that $\Omega^o(\mi,\mw,\me)$ is a transverse intersection at $(\mv,\me)$. Such points clearly form an open subset of $\muu(\mi,\mw)$. Now apply Lemma ~\ref{basic11} to conclude the proof.
\end{proof}
\section{Proof of the shift properties from Section ~\ref{termin}}\label{PP}
We are going to use the transversality result (Proposition ~\ref{carol}) in this section.
\begin{lemma}\label{Zariski}
Let $f:X\to Y$ be a quasi-projective morphism from an irreducible variety $X$ to a normal irreducible variety $Y$ of the same dimension. Let $y\in Y$ and $x_1,$ $\dots,$ $x_m$ isolated points of $f^{-1}(y)$. Then  $\deg(f)\geq m$.
\end{lemma}
\begin{proof}
Embed $X$ as an open subset of a variety $\bar{X}$ which is projective over $Y$:$$\xymatrix{X\ar@{^{(}->}[r]^{i}\ar[dr]^{f}&  \bar{X}\ar[d]^{\bar{f}}\\
                                           &     Y\\}$$
Let  $\bar{X}\leto{a} \tilde{Y}\leto{b} Y$ be the Stein factorization (~\cite{hart}, Section II.11.5) of
$\bar{f}$ where  $a$ has connected fibers and $b$ is finite. Therefore the images $a(x_1),\dots, a(x_s)$ are
distinct points of $\tilde{Y}$. Now use (~\cite{shafarevich}, Section II.6.3, Theorem 3) to conclude that
$\deg(f)\geq \deg(b)\geq m$.
\end{proof}

\begin{corollary}\label{Zariski2}\label{threepointtwo}
Let $\mi=(d,r,D,n,I)$ be a Schubert state with $\dimii=0$. Let $\mw=\mz_{D,n}$, and $\mf\in\Fl(\mw)$. Suppose that  there are $m$ isolated points $\mv_1,\dots,\mv_m$ in $\Omega^o(\mi,\mw,\mf)$. Then,  $\braI \geq m$.
\end{corollary}
\begin{proof}
It is easy to see that $\braI$ is the degree of the morphism (between smooth schemes) $\muu(\mi,\mw)\to \Fl(\mw)$. We can therefore use Lemma ~\ref{Zariski}.
\end{proof}
\subsection{Proof of Lemma ~\ref{jkl}}
The first claim follows from the equation $\cee(d,r,D,n)=\cee(d+r,r,D+n,n)$.

Clearly $\mz_{D,n}\tensor\mathcal{O}(-1)$ is isomorphic to $\mz_{D+n,n}$.
We therefore get an isomorphism between $\mfm(d,r,D,n)$ and $\mfm(r,d+r,D+n,n)$. There is also an natural isomorphism  $T:\Fl(\mz_{D,n})\leto{\sim}\Fl(\mz_{D,n}(-1))$. For $\mf$ in $\Fl(\mz_{D,n},\mpp)$, we find
a scheme theoretic isomorphism:
\begin{equation}\label{langl}
\Omega^o(\mi,\mz_{D,n},T(\mf))\leto{\sim} \Omega^o(\mj,\mz_{D,n},T(\mf)).
\end{equation}
 Therefore  Lemma ~\ref{Zariski2} implies that $\mi$ is not null $\Rightarrow$ $\mj$ is not null. The reverse implication is proved by tensoring with $\mathcal{O}_{\pone}(1)$.

In (3), if $\langle\mi\rangle >0$, pick a $\mf\in \Fl(\mz_{D,n})$  such that the set $\Omega^o(\mi,\mz_{D,n},\mf)$ has exactly $\langle\mi\rangle$ (reduced) points. Therefore  by Lemma ~\ref{Zariski2} and Equation ~\ref{langl}, $\langle\mj\rangle >\langle\mi\rangle $. Tensoring with $\mathcal{O}_{\pone}(-1)$ proves the reverse inequality.
\subsection{Proof of Proposition ~\ref{shiftt1}}
For $(1)$, we consider two cases
\begin{enumerate}
\item[(a)] $i_1=1$: In this case $$\codim(\omega_{J^{p}})=\codim(\omega_{I^{p}})- ({\emma}),\ \td=d-1$$ and hence $$\dim\mfm(d-1,r,D-1,n)=\dim\mfm(d,r,D,n)-({\emma}).$$ One verifies immediately that $\dim\mi\ =\ \dim S(p)\mi)$.
\item[(b)] $i_1\neq 1$: In this case $\codim(\omega_{J^{p}})=\codim(\omega_{I^{p}})+r$, $\td=d-1$ and $\dim\mfm(d,r,D-1,n)=\dim\mfm(d,r,D,n)+r$ and $\dim\mi=\dim S(p)\mi$.
 \end{enumerate}
 We are going to show that if $\mi$ is not null, then $S(p)\mi$ is not null and that $\langle\mi\rangle\leq \langle S(p)(\mi)\rangle.$ Since  $S(p)^n\mi=(d-r,r,D-n,n,I)$ where $\mi=(d,r,D,n,I)$, the use of iteration and Lemma ~\ref{jkl} finishes the proof of $(2)$ and $(3)$.

If $\mi$ is null then there is nothing to show. So assume that $\mi$ is not null. Let $\mw=\mz_{D,n}$ and choose a generic point $\mf=\prod_{q\in\mpp} F^q_{\bull}\in \Fl(\mw)$. We therefore have that $\Omega^o(\mi,\mw,\mf)$ is nonempty and of dimension $\dimm(\mi)$.

Choose a  uniformising parameter $t$ at $p$ and  define $\tilde{\mathcal{W}}$ to be the bundle whose sections are meromorphic sections $s$ of $\mw$ so that $ts$ is holomorphic and has fiber at $p$ in $F^p_1$ (this is a special case of the shift operation from Section ~\ref{shiftop}). Let $\mj=S(p)\mi$.

We claim $\tilde{\mathcal{W}}$ is isomorphic to $\mz_{D-1,n}$. The degree computation is clear. Its genericity has to be checked. Without loss of generality assume (tensor $\mw$ by an appropriate $\mathcal{O}_{\Bbb{P}^1}(l)$)
:
$$\mw=\mathcal{O}_{\Bbb{P}^1}^u\oplus\mathcal{O}_{\Bbb{P}^1}(-1)^{n-u}\  u\neq n.$$
The sub $\mathcal{S}=\mathcal{O}_{\Bbb{P}^1}^u$ is the part generated by global sections. Hence
we may assume (since our flags on $\mz_{D,n}$ are generic) that
$F^p_1\cap \ms_p=0$.

Let $\tilde{\mw}=\oplus_{l=1}^n\mathcal{O}(l)$. Since $\mw\hookrightarrow \tilde{\mw}$, we should have $a_l\geq -1$ for $l=1,\dots,r$. We claim $a_l\leq 0$ for $l=1,\dots,r$ and this would prove that $\tilde{\mw}$ is ES.
Assume the contrary, and let $\mathcal{O}_{\Bbb{P}^1}(1)\to \tilde{\mathcal{W}}$, be a nonzero map. We therefore find  a map $\mathcal{O}_{\Bbb{P}^1}\to \tilde{\mathcal{W}}(-p),$ or a section $s$ of $\tilde{\mathcal{W}}$ that vanishes at $p$.  Hence $s$ is a section of $\mathcal{W}$ whose fiber at $p$ is in $F^p_1$. But this is in contradiction to the assumption that there are no global sections of $\mw$ with fiber at $p$ in $F^p_1$.

We also induce complete flags on the fibers of $\tilde{\mw}$ at $\{p_1,\dots,p_s\}$ from those of $\mw$ (~\cite{b1}, Appendix) to obtain $\mg\in\Fl(\tmw)$. It is a standard fact that if we are given two vector bundles on a curve along with an isomorphism on a Zariski open set, then we obtain a set theoretic bijection between subbundles. This is via the saturation operation (see Section ~\ref{satu}). This is applied to $\mw$ and $\tilde{\mw}$ (and track is kept on the Schubert position of the fibers at $p$ as in ~\cite{b1}, Appendix)). We get that if $\mv\in\Omega^o(\mi,\mv,\mf)$, then $\tilde{\mv}$ obtained by  taking the saturation of  $\mv$ in $\tilde{\mw}$ is in $\Omega^o(\mj,\tilde{\mw},\mg)$ (and vice-versa).  Hence by Lemma ~\ref{threepointtwo}, $\mj$ is not null.

 The number of points in $\Omega^o(\mi,\mw,\mf)$ (for $\mf$ generic) is exactly $\langle\mi\rangle$ and therefore by Lemma ~\ref{Zariski2}, we have $\langle \mi\rangle\leq \langle\mj\rangle$.

 \appendix\section{Results from commutative algebra and representability criteria}
\subsection{Local criteria for flatness}
For the sake of reference, we collect some standard statements on flatness. The following lemma  is ~\cite{eisen}, Theorem 6.8, page 168.
\begin{lemma}\label{eisen}
Suppose that $(A,m)$ is a local Noetherian ring, and let $(S,n)$ be a local
 Noetherian $R$-algebra such that $mS\subseteq n$. If $M$ is a finitely generated $S$-module, then $M$ is flat as an $A$ module if and only if $\operatorname{Tor}^R_1(M,A/m)=0.$
\end{lemma}
The following lemma on flatness is ~\cite{akl}, Proposition V.3.4, page 94.
\begin{lemma}\label{sturm}
Let $R\to A$ and $A\to B$ be local homomorphisms of local rings and let $M$ be a finite $B$ module. Suppose that $A$ is flat over $R$. Then $M$ is flat over $A$ if and only if the following two conditions hold:
\begin{enumerate}
\item[(a)] $M$ is flat over $R$.
\item[(b)] $M\tensor_R k$ is flat over $A\tensor_R k$ where $k=R/m$ and $m$ is the maximal ideal.
\end{enumerate}
\end{lemma}
We also recall the standard fact: A finitely generated module over a Noetherian ring is flat if and only if it is locally free.
\begin{lemma}\label{mankl} Suppose that $T$ is a scheme, $\mw$ a vector bundle on $\pone$ and $\mv$ a coherent subsheaf of $\mw_T$. The following are equivalent:
\begin{enumerate}
\item  The quotient $\mq=\mw_T/\mv$ is flat over $T$.
\item  For any $t\in T$ (closed point),
$\mv_t\ \to\  \mw$ is an injection.
\end{enumerate}
Suppose the conditions above hold. Then, $\mv$ is a vector bundle on $\pone\times\  T$. Also, the following conditions are then equivalent:
\begin{enumerate}
\item[(a)] $\mq_t$ is a vector bundle for all $t\in T$
\item[(b)]$\mq$ is a vector bundle.
\item[(c)] $\mv$ is a subbundle (locally a direct summand) of $\mw_T$.
\end{enumerate}
\end{lemma}
\begin{proof}
Assume that $T=\Spec(A)$ where $A$ is local with maximal ideal $m$ and residue field $k=\kappa$. Let $t\in T$ correspond to the maximal ideal.
Consider the exact sequence
$$0\ \to\ \mv\ \to\ \mw_T\ \to\ \mq\ \to\  0.$$
Tensor this with $k(t)$ and obtain
$$0={\Tor}^1_A(\mw,k)\ \to\  {\Tor}^1_A(\mq_t,k)\ \to\  \mv_t\ \to\ \mw_t\ \to\ \mq_t\ \to\  0.$$
Therefore (2)$\Leftrightarrow$ (1) is clear by Lemma ~\ref{eisen}. Assuming that (1) and (2) hold, we see that $\mv$ is flat over $T$ (${\Tor}^2_A(\mq_t,k)\ \to\  {\Tor}^1_A(\mv_t,k)\ \to\  {\Tor}^1_A(\mw,k)$ is exact). $\mv_t$ is torsion free and hence  locally free on $\pone$. By Lemma ~\ref{sturm}, $\mv$ is flat over $\pone\times\  T$ and hence locally free.

For the equivalences $(a)\Leftrightarrow(b)\Leftrightarrow(c)$, the only non trivial implication is $(a)\Rightarrow(b)$, which follows from Lemma ~\ref{sturm}
and the flatness of $\mq$ over $T$.
\end{proof}
\begin{lemma}\label{injective}\label{aa}
Let $X$ be an irreducible variety, $S$ a scheme and  $\mv\ \to\  \mq$ a morphism of vector bundles on $X\times\ S$. Then, the subset of points $s\in S$ for which $\mv_s\ \to\ \mq_s$ is injective as a morphism of sheaves is open.
\end{lemma}
\begin{proof}
Taking exterior products one is immediately reduced to the case $\mv=\mathcal{O}$. We need to show that if $h$ is a global section of $\mq$ on $X\times\ S$, the subset $Z$ of $S$ formed by $s$ such that $h_s\neq 0$ in $H^0(X\times\  \{s\},\mq_s)$  is open.

But this is clear since if $h$ does not vanish at $(x,s)$, then there is a neighborhood $U$ of $s$ such that if $s'\in U$ then $h$ does not vanish at  $(x,s')$.
\end{proof}
\begin{lemma}\label{jacka}\label{barness}
Let $X$ and $S$ be schemes, with $X$ irreducible,  $\mv$ and $\mq$  vector bundles on $X\times\  S$ and $\phi:\mv\ \to\ \mq$ a morphism. The following are equivalent:
\begin{enumerate}
\item $\phi$ is injective and the cokernel of $\phi$ is flat over $S$.
\item For each closed point $s$ of $S$, $\phi_s$ is an injective as a morphism of sheaves on $\pone$.
\end{enumerate}
Furthermore, if $s\in S$, $\phi_s$ is injective (as a morphism of sheaves) iff it is injective on some nonempty open subset of $X$.
\end{lemma}
\begin{proof}
(1)$\Rightarrow$(2) is obvious from the definition of flatness. For
(2)$\Rightarrow$(1) we reason as follows: Break up $\phi$ into $2$ exact sequences $0\ \to\ \ms\ \to\ \mv\ \to\ \mt\ \to\  0$ and $0\ \to\ \mt\ \to\ \mq\ \to\ \mw\ \to\  0$. We have for $s\in
S$ a surjection $\mv_s\twoheadrightarrow \mt_s$ and a morphism $\mt_s\ \to\  \mq_s$ such that the composite is $\phi_s$. Therefore if $\phi_s$ is injective, $\mt_s\ \to\ \mq_s$ is
an injection and hence $\mw$ is flat over $S$. This implies (Lemma ~\ref{mankl}) that $\mt$ is a vector bundle and that $\ms$ is a vector bundle. But now we see that the rank of $\ms$ equals zero, so $\phi$ is injective.

The last statement is clear because $X$ is irreducible.
\end{proof}
\subsection{Some representability statements}\label{patrick}
Let $S$ be a scheme and let $E$ and $F$ be free $\mathcal{O}_S$-modules of rank $e$ and $f$, respectively. Let $\phi:E\to F$ be a morphism. The rank of $\phi$ is said to be $r$ if the cokernel is a locally free sheaf of rank $f-r$. With $\phi$ as above, if the cokernel of $\phi$ is locally free, then so are the image and the kernel.

With $S$, $E$, $F$ and $\phi$ as above, the $r$th degeneracy locus $D_r(\phi)$ of $\phi$ is the closed subset of all points $s\in S$ such that $\phi\tensor k(s)$ has rank less than or equal to $r$. We put a structure of a closed subscheme on $D_r(\phi)$ by writing it as the scheme of zeroes of the morphism
$$\bigwedge^{r+1}E\to \bigwedge^{r+1}F.$$
Notice that the subschemes $Z_r(\phi)=D_r(\phi)-D_{r-1}(\phi)$ partition $S$ into a disjoint union of locally closed subschemes.

\begin{lemma} With notation as above., if $T=Z_r(\phi)$, the map $\pi_T:E_T\to F_T$, has rank $r$. That is, the cokernel is locally free of rank $f-r$.
\end{lemma}
\begin{proof}See ~\cite{eisen}, Corollary 20.5, Proposition 20.8.
\end{proof}
Given $S$, $E$, $F$ we form the scheme $\Hom(E,F)$ over $S$. Write $\pi:\Hom(E,F)\to S$ for the structure map. There is a natural map $\phi:\pi^{*}E\to \pi^{*}F$. we write $\Hom_r(E,F)$ for the locally closed subscheme $Z_r(\phi)$. This represents a functor:
\begin{lemma}
Let S be a Noetherian scheme; then $\Hom_r(E,F)$ represents the functor:
$$T\leadsto \{\psi:E_T\to F_T\mid\rk(\psi)=r\}$$
\end{lemma}
\begin{proof}Standard.
\end{proof}

Let $Z$ be a scheme with a vector bundle $\mv$ of rank $r$ which is filtered by a series of subbundles:
$$\mv^{1}\subseteq \dots \subseteq \mv^{k+1}=\mv$$
where $\mv^{l}$ is of rank $b_l$ for $l=1,\dots,k+1$.
 Let $a_0=0\leq a_1\leq\dots\leq a_{k+1}= r$ be nonnegative integers. Consider the functor $G$ whose value $G(T)$ over a  scheme $T$ over $Z$ is the set  of  complete filtrations by subbundles:
$$0\subsetneq\mw_1\subsetneq\dots\subsetneq \mw_r=\mv_T$$
so that $\mv_T^{l} \subseteq \mw_{a_l}$ for $l=1,\dots,k$.
\begin{lemma}\label{basic1}
The functor $G$ is representable by the topmost element in  a tower of Grassmann bundles over $Z$. The fiber dimension of the representing scheme over $Z$ is $$\frac{1}{2}r(r-1)-\sum_{l=1}^k(a_{l+1}-a_l)b_l.$$
\end{lemma}
\begin{proof}
Let  $X_k=\Gr(a_k-b_k,\frac{\mv}{\mv^{k}})$. On $X_k$ there is a natural bundle $\mw_{a_k}$ of rank $a_k$ which contains $p_k^*(\mv^{k})$. Now form the Grassmann bundle $X_{k-1}=\Gr(a_{k-1}-b_{k-1},\frac{\mw_{a_k}}{p_k^*\mv^{k-1}})$ over $X_k$ and iterate this procedure producing a tower of Grassmann bundles:
$$X_1\leto{p_1}X_2\leto{p_2}\dots\leto{p_{k-1}}X_k\leto{p_k} Z$$
Let $q:X_1\to Z$.
Clearly,
$$\dim(X_1)=\sum_{l=1}^{k}(a_l-b_l)(a_{l+1}-a_l).$$
We have gaps in the filtration of $q^*\mv$ (not every $\mw_i$ appears as a $\mw_{a_l}$), which we ``fill'' arbitrarily by forming suitable flag bundles. Therefore the representing scheme is of dimension (over $Z$)
$$\sum_{l=1}^{k}[(a_l-b_l)(a_{l+1}-a_l)+\dim(\operatorname{Fl}(\kappa^{a_{l+1}-a_{l}}))]+\frac{a_1(a_1-1)}{2}$$
$$=\{\sum_{\ell=0}^k[a_{\ell}+ (a_{\ell} +1)+\dots+(a_{\ell}+a_{\ell+1}-a_{\ell}-1)]\} -\sum_{l=1}^r(a_{l+1}-a_l)b_l.$$

The first bracketed term is $1+2+\dots +(a_{k+1}-1)=\frac{r(r-1)}{2}$and the proof is complete.
\end{proof}
\subsection{More representability}
Let $X$ and $S$ be varieties, $p:X\to S$ a morphism and $\mv$ a coherent sheaf on $X$ which is flat over $S$. Let $F$ be the contravariant functor: schemes/$S$ to abelian groups, which assigns to a scheme $T$ over $S$, the abelian group $F(T)=H^0(X_T,\mv_T)$ (where $\mv_T$ is the pull back of $\mv$ to $X_T=X\times_S T$ via the morphism $X_T \to X$).
\begin{lemma}\label{somelemma}
In the situation above assume that
\begin{enumerate}
\item[(a)] $p_* \mv$ is a vector bundle  on $S$.
\item[(b)] For any morphism $f:T\to S$, if
\begin{equation}\label{somenumber}
\xymatrix{X_T\ar[d]^{p'}\ar[r]^{f'} &  X\ar[d]^{p}\\
            T\ar[r]^{f}               &  S}
\end{equation}
\end{enumerate}
 is the ``base change'' cartesian square, the natural morphism $f^*p_*\mv\to p'_*f'^*\mv$ is an isomorphism.

Then the functor $F$ above is represented by the total space of the vector bundle $p_*\mv$.
\end{lemma}
\begin{proof} We first note that if $\mw$ is a vector bundle on a scheme $Y$
with total space $W$, then $\homo_Y(Y,W)=H^0(Y,\mw).$

Now, let $\ma=p_*\mv$. Let $A$ be the total space of this vector bundle. Let $f:T\to S$ be a morphism and use notation from the Cartesian square ~\ref{somenumber}. We have $F(T)=H^0(X_T,f'^*\mv)=H^0(T,p'_{*}f'^{*}\mv)$. But the latter is (naturally) isomorphic to $H^0(T,f^*\ma)$ by the second assumption. Therefore if $A_T$  is the total space of $f^*\ma$, we have $F(T)=\homo_T(T,A_T)$. The lemma follows from the following cartesian square:
$$\xymatrix{A_T\ar[r]\ar[d] & A\ar[d]\\
              T\ar[r] &             S}$$
Hence $F(T)$ is identified with the set of morphisms $T\to A$ (over $S$).
\end{proof}
Conditions (a) and (b) in Lemma ~\ref{somelemma} hold  in each of the following two cases:
\begin{enumerate}
\item $S=\Spec(\kappa)$. In this case the second hypotheses is true because ``cohomology commutes with flat base change'' (see ~\cite{hart}, Chapter III, Proposition 9.3).
\item The higher direct images of $\mv$ vanish. In fact if we assume that for all points $s\in S$, $\Hone(p^{-1}(s),\mv_s)=0$, then the hypothesis are valid (cf. ~\cite{curves}, Corollary 1 from Lecture 7).
\end{enumerate}
\begin{corollary}\label{rep1}
Let $\mv$, $\mq$ be vector bundles on a scheme $X$. Then, the vector space $Hom(\mv,\mq)$ considered as a scheme over $\Spec(\kappa)$ represents the functor: schemes/$\kappa$ to (sets) given by
 $T\leadsto \homo_{X\times T}(\mv_T,\mq_T)$.
\end{corollary}
\subsection{Shift operations}\label{shiftop}
Fix a point $p\in \pone$. We have a natural inclusion of sheaves $\mathcal{O}_{\pone}\subseteq \mathcal{O}_{\pone}(p)$ and given a uniformising parameter $t$ at $p$, we have in a neighborhood of $p$, an isomorphism $\mathcal{O}_{\pone}(p)\leto{t}\mathcal{O}_{\pone}$.

Let $\mv$ be a vector bundle on $\pone\times\ Y$, $p\in \pone$ and $\mb\subseteq \mv_p$ a subbundle ($\mv_p$ is a bundle on $p\times\  Y$). There is a canonical inclusion $\mv\subseteq\mv(p)=\mv\tensor p_Y^*\mathcal{O}(p)$.

We define a new sheaf $\tmv$ such that $\mv\subseteq\tmv\subseteq\mv(p)$: $\tmv$ coincides with $\mv$ outside of $(\pone-p)\times\  Y$ and in a neighborhood of $p\times Y$, is the kernel of the composite: $ \mv(p)\leto{t} \mv\to (\mv_p/\mb)$ where the first map is multiplication by $t$. Clearly, $\tmv$ is  independent of the choice of the uniformising parameter $t$.

It is useful to have a local model of this operation: If  $y\in Y$, there is an open subset $U_y$ of $Y$
containing $y$, a neighborhood $\tilde{U}$ in $X\times\ Y$ of $p\times\ U_y$, and a  splitting $\mv=\mc\oplus\mt$ on
$\tilde{U}$ such that $\mc$ restricted to $p\times\  U_y$ is the same as $\mb$. On $\tilde{U}$,   $\tmv$ is
identified with the subsheaf $\mc(p)\oplus\mt$ of $\mv(p)$. Therefore $\mv\to \tmv$ has cokernel $\mb(p)$, and
$\mv_p\to\tmv_p$ has kernel $\mb$.

We can perform the shift operation at various points of a bundle simultaneously. Let $\mv$ be a vector bundle on $\pone\times\ Y$ and for  $p\in\mpp$, suppose $\mb=\prod_p \mb^p$ where $\mb^p \subseteq \mv_p$ is a subbundle. Then we form the shift $\mv\subset\Psi(\mv,\mb)\subset\mv\tensor\mathcal{O}(\sum_{p\in \yess}p)$.
\begin{lemma}
In the above setup with $\mb=\prod_{p\in \mpp}\mb^p$, and $\tmv=\Psi(\mv,\mb)$,
\begin{enumerate}
\item $\tmv$ is locally free, the quotients $\tmv/\mv$ and $\mv(p)/\tmv$ are flat over $Y$, and  the formation of $\tmv$ commutes with base change in $Y$.
\item For $p\in\mpp$, $\mv\to \tmv$ has cokernel $\mb^p\otimes\mathcal{O}(p)$, and $\mv_p\to\tmv_p$ has kernel $\mb^p$.
\item For any vector bundle $\mq$ on $\pone\times\   Y$, $Hom(\tmv,\mq)\hookrightarrow\Hom(\mv,\mq)$ and there is a 1-1 correspondence between the following objects:
\begin{enumerate}
\item Morphisms $\phi:\mv\to\mq$ such that $\mb^p$ is in the kernel of the map $\phi_p:\mv_p\to \mq_p$ for each $p\in\mpp$,
\item Morphisms  $\phi': \tmv\to \mq$
\end{enumerate}
Moreover, for $\phi$ as in (a), $\phi$ is injective with a flat cokernel if and only if  the same is true for $\phi'$.
\end{enumerate}
\end{lemma}
\begin{proof}
For ease of exposition assume $\mpp=\{p\}$ and suppress the superscript $p$ in the notation. Let $t$ be a uniformising parameter at $p$. From the local model, it is easy to see $(1)$ and  $(2)$ hold. We have an exact sequence of the form
$$0\to\mv\to\tmv\to\mb(p)(=\mb\tensor \mathcal{O}(p))\to 0$$
and hence an exact sequence
$$0\to\homo(\mb(p),\mq)\to\homo(\tmv,\mq)\to\homo(\mv,\mq).$$
The first term is clearly $0$ since multiplication by $t$ is injective on the second and third terms. We have therefore proved the injection part of (3).

In (3), it is clear that given $\phi'$ as in (b), we can obtain $\phi$ by composition $\mv\to\tmv\to Q$. The kernel of $\phi_p$ contains the kernel of $\mv_p\to \tmv_p$ which is $\mb$ .

For the other direction, let $\phi$ be a morphism as in $3(a)$ and consider the following diagram,
$$\xymatrix{&   &     \tmv\ar[d]\\
           & \mv\ar[r]\ar[d]^{\phi}\ar[ur] &  \mv(p)\ar[r]^{t}\ar[d]^{\phi(p)} & \mv_p/\mb\ar[d]^{\phi_p} \\
          0\ar[r] &  \mq\ar[r]       &\mq(p)\ar[r]^{t} & \mq_p\ar[r] & 0.\\
}$$

 Note that $\phi_p:\mv_p\to \mq_p$ factors through $\mv_p/\mb$ because of the hypothesis in (a). The bottom row is an exact sequence and the composite $\tmv\to \mq_p$ is zero because it factors through the $0$ map $\tmv\to\mv_p/\mb$. We therefore find a factorization  $\phi':\tmv\to \mq$. The last statement follows from Lemma ~\ref{jacka}.

\end{proof}
\subsection{Parabolic bundles}\label{satu}
We consider parabolic bundles on $\pone$ with parabolic structure at the points $\mpp$.

A parabolic bundle $\underline{\mw}=(\mw,\me,w)$ on $(\pone,\mpp)$ is a $(D,n)$-vector bundle $\mw$ on $\pone$, a collection of complete flags $\me=\prod_{p\in\mpp}E^p_{\bull}$ and a function
$$w:\mpp\times\{1,\dots,n\}\to\Bbb{R}$$
such that, denoting $w(p,a)$ by $w^p_a$, we for each $p\in\mpp$,
$$w^p_1\geq w^p_2\geq\dots\geq w^p_n\geq w^p_1-1.$$
We define $w^p_0=w^p_n+1$.

 For a parabolic bundle $\umw$ as above and a subbundle $\mv\subseteq\mw$, let $\mi=(d,r,D,n,I)$ be the Schubert state determined from $\mv\in\Omega^o(\mi,\mw,\me)$.
\begin{enumerate}
\item Define the weight of $\mv$ by
$$\wt(\mv,\umw)=\sum_{p\in\mpp}\sum_{a\in I^{p}}w^p_a$$
\item The parabolic weight of $E$ is defined to be
$$\pa(\mv,\umw)= -d +\wt(\mv,\umw).$$
\item The parabolic slope of $\mv$ is defined to be
$$\mupar(\mv,\umw)=\frac{\pa(\mv,\umw)}{r}$$
\end{enumerate}
Let $\mw$ be a vector bundle on a smooth curve $C$ and $\mv\subseteq\mw$ a coherent subsheaf. Let $T$ be the torsion subsheaf of $\mw/\mv$ and $\tmv$ its inverse image in $\mw$. $\tmv$ is a subbundle of $\mw$ containing $\mv$ and is called the {\em saturation} of $\mv$ in $\mw$.
\begin{lemma}\label{naive} Let $\umw$ be a parabolic bundle as above and $\phi:\mv\hookrightarrow\mw$ a coherent subsheaf of degree $-d$ and rank $r$ with saturation $\tmv$. Suppose that there is a  $\mg=\prod_p G^p_{\bull}\in\Fl(\mv)$ and a function
$$\gamma:\mpp\times\{1,\dots,r\}\to \{0,1,\dots,n\}$$
such that
$$\phi_p(G^p_{a})\subseteq E^p_{{\gamma}^p(a)},\text{ }a=1,\dots,r, p\in\mpp$$
Then
$$\pa(\tmv,\umw)\geq -d +\sum_{p\in\mpp}\sum_{a=1}^r w^p_{{\gamma}^p_a}.$$
\end{lemma}
\begin{proof}
Let $\mb(p)$ be the kernel of $\phi_p$ and $b(p)$ the rank of $\mb(p)$. Let $\mi=(\td,r,D,n,I)$ be the Schubert state determined from $\tmv\in \Omega^o(\mi,\mw,\me)$. It is easy to see that $\deg(\mv)\geq -d+\sum_{p\in\mpp}b(p)$ and hence
$$\pa(\tmv,\umw)\geq -d+\sum_p b(p)+\sum_{p\in\mpp}\sum_{a\in I^{p}}w^p_a.$$
Lemma ~\ref{naive} now follows from:
\begin{claim}
For each $p\in\mpp$,
$$b(p)+\sum_{a\in I^{p}}w^p_a \geq \sum_{a=1}^rw^p_{{\gamma}^p(a)}.$$
\end{claim}
For the claim, we fix a $p$ and assume $\mb(p)\in \Omega^o_{H}(G^p_{\bull})$. Let $H=\{h_1<\dots<h_{b(p)}\}$ and $\{u_1<\dots<u_{r-b(p)}\}=\{1,\dots,r\}-H.$ Let $I^p=\{i_1<\dots<i_r\}$. We see that
$\phi_p(G^p_{u_t})$ is at least $t$ dimensional so $i^p_t\leq {\gamma}^p(u_t)$ for $t=1,\dots, r-b(p)$.
$$\sum_{a\in I^{p}}w^p_a - \sum_{a=1}^rw^p_{{\gamma}^p(a)}\geq \sum_{a=r-b(p)+1}^r w^p_{i_a} - \sum_{{\ell}=1}^{b(p)} w^p_{{\gamma}^p({h_{\ell}})}$$
$$=\sum_{{\ell}=1}^{b(p)} [w^p_{i_{r-{b(p)}+{\ell}}} - w^p_{{\gamma}^p({h_{\ell}})}]\geq\sum_{{\ell}=1}^{b(p)} (-1)= -{b(p)}.$$
\end{proof}

\bibliographystyle{plain}
\def\noopsort#1{}

\end{document}